\documentclass[11pt]{article}
\usepackage{latexsym}
\usepackage[all]{xy}
\usepackage{amsmath}
\usepackage{amssymb}
\usepackage{amsfonts}
\usepackage[applemac]{inputenc}
\RequirePackage{hyperref}

\newtheorem{thm}{Theorem}[section]
\newtheorem{prop}[thm]{Proposition}
\newtheorem{lem}[thm]{Lemma}

\newtheorem{df}[thm]{Definition}
\newtheorem{cor}[thm]{Corollary}

\newtheorem{rmk}[thm]{Remark}
\newtheorem{ex}[thm]{Example}
\newtheorem{exs}[thm]{Examples}

\newcommand{\car}{\mathrm{char}}

\newcommand{\lperf}{\mathsf{L}_{\mathsf{Perf}}}
\newcommand{\lqcoh}{\mathsf{L}_{\mathsf{QCoh}}}

\newcommand{\dgmod}{\mathsf{dgmod}}
\newcommand{\cdga}{\mathsf{cdga}^{\leq 0}_{k}}
\newcommand{\dga}{\mathsf{dga}}
\newcommand{\salg}{\mathsf{salg}_{k}}
\newcommand{\dgam}{\mathsf{dga}^{\leq 0}_{k}}
\newcommand{\com}{\mathsf{C}^{\leq 0}(k)}
\newcommand{\co}{\mathsf{C}(k)}

\newcommand{\idgmod}{\underline{\mathsf{dgmod}}}
\newcommand{\icdga}{\underline{\mathsf{cdga}}^{\leq 0}_{k}}
\newcommand{\idga}{\underline{\mathsf{dga}}}
\newcommand{\isalg}{\underline{\mathsf{salg}}_{k}}
\newcommand{\idgam}{\underline{\mathsf{dga}}^{\leq 0}_{k}}
\newcommand{\icom}{\underline{\mathsf{C}}^{\leq 0}(k)}
\newcommand{\ico}{\underline{\mathsf{C}}(k)}
\newcommand{\iperf}{\underline{\mathsf{Perf}}}

\newcommand{\ssets}{\mathsf{SSets}}
\newcommand{\dst}{\mathsf{dSt}_{k}}
\newcommand{\spec}{\mathbb{R}\mathrm{Spec}}

\newcommand{\ho}{\mathrm{Ho}}

\begin{document}

\title{\textbf{Quadratic forms and Clifford algebras \\ on derived stacks}}
\bigskip
\bigskip
\bigskip

\author{
\bigskip \\ Gabriele Vezzosi\\
\small{Dipartimento di Matematica ed Informatica} \\ 
\small{Firenze - Italy}
\\}


\maketitle

\begin{abstract}
\noindent In this paper we present an approach to quadratic structures in derived algebraic geometry. We define derived $n$-shifted quadratic complexes, over derived affine stacks and over general derived stacks, and give several examples of those. We define the associated notion of derived Clifford algebra, in all these contexts, and compare it with its classical version, when they both apply. Finally, we prove three main existence results for derived shifted quadratic forms over derived stacks, define a derived version of the Grothendieck-Witt group of a derived stack, and compare it to the classical one.
\end{abstract}

\noindent \textbf{Keywords -} Derived algebraic geometry, quadratic forms, Clifford algebras

\tableofcontents

\bigskip
\bigskip

\section*{Introduction} Prompted by the recent introduction of \emph{symplectic forms} in derived algebraic geometry (\cite{ptvv}), we present a global approach to the theory of \emph{quadratic forms} on very general moduli spaces, called derived stacks. In the case where the derived stack is just the spectrum of a ring $k$ where $2$ is invertible, these quadratic forms are defined on complexes $C$ of $k$-modules, and are maps, in the derived category of $k$, from  $\mathrm{Sym}_{k}^{2} C$ to $k[n]$, where $\mathrm{Sym}_{k}^{2}(-)$ denotes the derived functor of the second symmetric power over $k$. There is an obvious notion of non-degeneracy for such a quadratic form, saying that the induced adjoint map $C \to C^{\vee}[n]$ is an isomorphism in the derived category of $k$. \\ The derived features are therefore two: first of all the map is a morphism in the \emph{derived category} of $k$, and secondly, and most importantly, we allow for a \emph{shift} in the target. These features accommodate for various symmetric shifted duality situations in topology, the motivating one being classical Poincaré duality. \\We present a \emph{globalization} of the above particular case to quadratic forms on Modules over a derived stack that uses in an essential way the refined (i.e. homotopical or, equivalently, $\infty$-categorical) features of derived algebraic geometry. When the Module in question is the tangent complex, we obtain what we call \emph{shifted quadratic stacks}. We remark that the main definitions of derived quadratic forms and derived quadratic stacks are slight modifications of the notion of derived symplectic structure from \cite{ptvv} (without the complication coming from closedness data). In particular, we are able to reproduce in the quadratic case, two of the main existence theorems in \cite{ptvv}: the existence of a shifted quadratic form on the stack of maps from a $\mathcal{O}$-compact, $\mathcal{O}$-oriented derived stack to a shifted quadratic stack (Theorem \ref{map}), and the existence of a quadratic form on the homotopy fiber product of two null-mappings to a shifted quadratic stack (Theorem \ref{lagr}).  As a consequence, we get that the derived looping of a shifted quadratic stack decreases the shift by one (Corollary \ref{loop}). We also observe that any shifted symplectic structure gives rise to a shifted quadratic Module. \\ The third main topic developed in the paper, after the general theory of quadratic forms on derived moduli spaces and the main existence theorems, is the definition and study of a new derived version of the \emph{Clifford algebra} associated to a shifted quadratic Module over a derived stack (so, in particular, to any shifted quadratic stack). We prove various basic properties of this derived Clifford algebra and give a theorem comparing it to the classical Clifford algebra, when they are both defined: the classical Clifford algebra happens to be the truncation at $\mathrm{H}^0$ of the derived one (see Corollary \ref{comparazione}). We also introduce the notion a shifted derived version of the Grothendieck-Witt group of an aribitrary derived stack $X$, and compare it to the classical counterpart when $X= \mathbb{R}\mathrm{Spec}\, k$ and the shift is $0$ (Proposition \ref{gwcomp}). \\
An appendix establishes the basics of the homotopy theory of (externally) $\mathbb{Z}/2$-graded differential graded algebras needed in the main text to properly define the $\mathbb{Z}/2$-grading on the derived Clifford algebra.\\		
Some leftovers and future directions of this work are discussed below, in this introduction.\\

\noindent \textbf{Alternative approaches}. An alternative approach to derived Clifford algebras with respect to the one used in this paper is the following, originally due in the underived and unshifted case to A. Preygel (\cite{pre}). It applies only in the even-shifted case, exactly as the one presented in this paper. Let $(E, q)$ be a $2n$-shifted quadratic complex on $X$ (see Definition \ref{dq} and \ref{global}). Here $X$ can be any derived Artin stack locally of finite presentation over the base ring $k$, but the reader could stick to the derived affine case of Def. \ref{dq} without loosing any essential feature. The derived quadratic form $q$ induces on the linear derived stack $$\mathbb{V}(E[-n]):= \mathbb{R}\mathbf{Spec} (\mathrm{Sym}_{\mathcal{O}_{X}}(E^{\vee}[n]))$$ a global function $f_{q}: \mathbb{V}(E[-n]) \to \mathbb{A}_{k}^1$. Then we may view the pair $(\mathbb{V}(E[n]), f_{q})$ as a \emph{Landau-Ginzburg pair} (see \cite{pre}). One can then extract from the corresponding matrix factorization category $\mathbf{MF}((\mathbb{V}(E[n]), f_{q}))$, a derived Clifford algebra for $(E,q)$, via a (derived, shifted) variant of \cite[Thm. 9.3.4]{pre}. Its $\mathbb{Z}/2$-grading appears here in a natural way. Note that for derived quadratic stacks, i.e. when $E= \mathbb{T}_{X}$, then $\mathbb{V}(E[-n])$ is  exactly the $(-n)$-shifted derived  tangent stack of $X$.\\

\noindent \textbf{Leftovers}.  Many interesting topics are not present in this first paper on the derived theory of quadratic forms. Namely, a treatment of involutions, analogues of classical classification theorems (like Witt cancellation), more details and applications of the derived Grothendieck-Witt group (defined in Section \ref{GW}), and a more thorough investigation of the  relations to derived Azumaya algebras (\cite{deraz}), especially in the $\mathbb{Z}/2$-graded case where we may expect a map from the derived Grothendieck-Witt group  to the derived Brauer-Wall group (classifying  derived $\mathbb{Z}/2$-graded Azumaya up to $\mathbb{Z}/2$-graded Morita equivalences). We plan to come back at least to some of these leftovers in a future paper.\\

\noindent \textbf{Generalizations.} One natural, mild, extension of the theory developed in this paper can be obtained by working over a base ring $k$ \emph{with involution} (along the lines of Ranicki's work \cite{ra}). This would allow to treat, for example, derived (shifted) hermitian forms, objects that arise naturally (together with closed shifted differential forms) in derived K\"ahler geometry, a topic yet to be investigated.\\
A wider generalization of the present theory involves a categorification process: define a derived quadratic form on a dg-category. In order to do this one needs a homotopy invariant notion of the second symmetric power of a dg-category, i.e. a slight reformulation of the notion introduced by Ganter and Kapranov in \cite{gaka}. As further explained by M. Schlichting, the role of the shifted $k[n]$ will be taken here by $\mathcal{C}^{[n]}_{k}$ (see  \cite[\S 1.9]{sch}). We plan to come back to this categorification in a future paper. One could expect e.g. that the quasi-coherent dg-category of a derived quadratic stack carries a derived quadratic form, but it might be the case that the converse is not true. This would of course make the categorification more interesting.\\

\noindent \textbf{Acknowledgements.} The atmosphere, food and waves in Portugal, gave a perfect start to this work.\\ On a more intellectual level, this paper grew out of some questions I raised during the preparation of the paper \cite{ptvv}, while I was visiting the IHES: I had interesting start-up conversations on the subject with Bertrand To\"en, and Maxim Kontsevich was the first who thought that a derived notion of the Clifford algebra could be a useful object. I thank both of them. One more thank to Bertrand To\"en for various useful comments on a first draft of the paper,  to Luca Migliorini for his explanations about intersection (co)homology,  to Nick Rozenblyum for helpful conversations related to the derived  Grothendieck-Wtt group, and to Tony Pantev for his interest in this work and his comments. I was also influenced by Andrew Ranicki's work on algebraic surgery, and by the notes of Jacob Lurie's Harvard Course on Algebraic L-theory and Surgery (Spring 2011).\\

\bigskip

\noindent \textbf{Notations and conventions.} 

\begin{itemize}
\item  $k$ will denote our base commutative ring such that $2\neq0$ in $k$. When $k$ will be needed to be a field, we will also use the alternative notation $k=\mathbb{F}$.
\item When we say that $(V,Q)$ is a (classical) \emph{quadratic module} over $k$, we actually mean that $V$ is a $k$-module, and that $Q$ is a symmetric bilinear form on $V$ over $k$ (these might bear the name of symmetric bilinear modules, but this name is almost never used in the literature). Note that classical quadratic modules, as defined e.g. in \cite{mire}, coincide with our if $2$ is invertible in $k$.
\item By $\infty$-category, we mean a Segal category (so in particular, any simplicial category will be viewed as an $\infty$-category), and we use the results and notations from \cite{chern}. Equivalently, one might work in the framework of quasi-categories (see \cite{htt}). For any $\infty$-category $T$, and any pair of objects $(x,y)$ in $T$, we denote by $\mathrm{Map}_{T}(x,y)$ the corresponding mapping space (well defined in $\ho (\ssets)$).
\item As a general rule, $\infty$ categories will be denoted in underlined fonts like $\underline{\mathsf{C}}$, and the corresponding underlying model categories (whenever they exist) with no underline, like $\mathsf{C}$: if $\mathsf{C}$ is a model category then the associated $\infty$-category $\underline{\mathbf{C}}$ is the Dwyer-Kan simplicial localization (\cite{dk}) of $\mathsf{C}$ with respect to weak equivalences (or its homotopy coherent nerve \cite[1.1.5]{htt}, if one prefers working within the context of quasi-categories).
\item $\co$ will denote the model category of unbounded cochain complexes of $k$-modules with surjections as fibrations, and quasi-isomorphisms as equivalences. It is a symmetric monoidal model category with the usual tensor product $\otimes_{k}$ of complexes over $k$, and it satisfies the monoid axiom (\cite{ss1}). The corresponding $\infty$-category will be denoted by $\ico$.
\item $\com$ will denote the model category of cochain complexes of $k$-modules concentrated in non-positive degrees, with surjections in strictly negative degrees as fibrations, and quasi-isomorphisms as equivalences. It is a symmetric monoidal model category with the usual tensor product $\otimes_{k}$ of complexes over $k$,  and it satisfies the monoid axiom (\cite{ss1}). The corresponding $\infty$-category will be denoted by $\icom$.
\item $\cdga$ denotes the category of differential \emph{non-positively} graded algebras over $k$, with differential increasing the degree by 1. 
If $\car k = 0$, then one might endow $\cdga$ with the usual model structure for which fibrations are surjections in negative degrees, and equivalences are quasi-isomorphisms (see \cite[\S 2.2.1]{hagII}). The associated $\infty$-category will be denoted by $\icdga$. In general, i.e. for $\car \, k $ not necessarily zero, we will just work with $\icdga$ defined as the simplicial localization of $\cdga$ along quasi-isomorphisms. 
\item $\dga_{k}$ denotes the category of unbounded differential graded cochain algebras over $k$ (not necessarily commutative). 
We will always consider $\dga_{k}$ endowed with the usual model structure for which fibrations are surjections, and equivalences are quasi-isomorphisms (see \cite{ss1}). Its associated $\infty$-category will be denoted by $\idga_{k}$.
\item $\dgam$ denotes the category of differential non-positively graded cochain algebras over $k$. 
We will always consider $\dgam$ endowed with the usual model structure for which fibrations are surjections in strictly negative degrees, and equivalences are quasi-isomorphisms (see \cite[\S 2.3]{hagII}). Its associated $\infty$-category will be denoted by $\idgam$.
\item $\salg$ denotes the model category of simplicial commutative $k$-algebras where weak equivalences and fibrations are detected on the underlying morphisms of simplicial sets (see \cite[\S 2.2]{hagII}). Its associated $\infty$-category will be denoted by $\isalg$.
\item We will denote by $\mathcal{S}$ the $\infty$-category of spaces or simplicial sets, i.e. $\mathcal{S}$ is the Dwyer-Kan localization of the category $\ssets$ with respect to weak equivalences.
\item $\dst$ will denote the model category of derived stacks over $k$ (\cite[\S 2.2]{hagII}). Its associated $\infty$-category will be denoted by $\underline{\dst}$.
\item A derived geometric stack over $k$ (see \cite[\S 2.2]{hagII}) will be called \emph{lfp} if it is locally of finite presentation over $k$. In particular, any lfp derived geometric stack $X \in \dst$ has a cotangent complex $\mathbb{L}_{X}$ that is perfect over $\mathcal{O}_X$.  
\end{itemize}

\section{Derived quadratic complexes}

If $A\in \salg$, we denote by $\mathrm{N}(A)$ its normalization (see \cite{dold}). Since The normalization functor is lax symmetric monoidal, $\mathrm{N}(A)\in \cdga$ (see \cite{ss2}) and we will write 
\begin{itemize}
\item $A-\dgmod$ for the model category of unbounded dg-modules over $\mathrm{N}(A)$. This is a symmetric monoidal model category satisfying the monoid axiom (\cite{ss1}). The corresponding $\infty$-category will be denoted as $A-\underline{\dgmod}$. We will write $A-\iperf$ for the full $\infty$-subcategory of perfect dg-modules over $\mathrm{N}(A)$.
\item $\mathsf{D}(A):= \ho (A-\dgmod)$ is the unbounded derived category of $\mathrm{N}(A)$, $\otimes_{A}^{\mathbb{L}}\equiv \otimes_{A}$ its induced monoidal structure, and $\mathbb{R}\underline{\mathrm{Hom}}_{A}$ the adjoint internal Hom-functor.
\item $A-\dga$ the model category of differential graded algebras over $A$, i.e. the category of monoids in the symmetric monoidal model category $A-\dgmod$ (see \cite{ss2}). The corresponding $\infty$-category will be denoted as $A-\idga$.
\item For any $n\in \mathbb{Z}$, $A[n]:= \mathrm{N}(A)[n]$ as an object in $A-\dgmod$.
\item $\bigwedge^2 \equiv \bigwedge^2_k$, $\otimes \equiv \otimes_{k}$, $\mathrm{Sym}_{k}^2\equiv \mathrm{Sym}_{k}^2$, $\bigwedge^2_A$, $\otimes_{A}$, $\mathrm{Sym}_{A}^2$ (for $A\in \salg$, identified with its normalization cdga $\mathrm{N}(A)$) will always be \emph{$\infty$-functors}. \\ In particular, $\mathrm{Sym}_{A}^2 : A-\underline{\dgmod} \longrightarrow A-\underline{\dgmod} $ is defined as follows\footnote{If $n!$ is invertible in $k$, a completely analogous definition holds for the $\infty$-functor $Sym_{A}^{n}$.}. First note that $A-\underline{\dgmod}$ is also equivalent to the Dwyer-Kan localization of cofibrant-fibrant objects in $A-\dgmod$ (and all objects in $A-\dgmod$ are fibrant), and $P \mapsto P\otimes_{\mathrm{N}(A)}P$ preserves cofibrant objects (by the monoid axiom) and weak equivalences between (since cofibrant objects are homotopically flat). So the assignment $^{\mathbb{L}}\otimes_{A}^2: P \mapsto QP\otimes_{\mathrm{N}(A)}QP$, where $Q$ denotes a cofibrant replacement functor in $A-\dgmod$ defines, after simplicial localization, an $\infty$-functor $^{\mathbb{L}}\otimes_{A}^2: A-\underline{\dgmod} \longrightarrow A-\underline{\dgmod}^{\Sigma_2}$ where $A-\underline{\dgmod}^{\Sigma_2}$ is the $\infty$-category obtained by Dwyer-Kan localizing the projective model structure of $\mathrm{N}(A)$-dg modules endowed with a $\Sigma_2$-action (i.e. the projective model structure on $B\Sigma_2$-diagrams in $A-\dgmod$). 
Finally, the $\Sigma_2$-coinvariants functor $(-)_{\Sigma_2}: A-\dgmod^{\Sigma_2} \longrightarrow A-\dgmod$ preserves weak equivalences since $2$ is invertible in $k$, and therefore the coinvariants functor from $k$-modules linear representations of $\Sigma_2$ to $k$-modules is exact. Thus, we obtain, after Dwyer-Kan localization, an $\infty$-functor $(-)_{\Sigma_2}: A-\underline{\dgmod}^{\Sigma_2} \longrightarrow A-\underline{\dgmod}$. Then, $\mathrm{Sym}_{A}^2$ is defined by $$\mathrm{Sym}_{A}^2 \simeq (-)_{\Sigma_2} \circ \,\,^{\mathbb{L}}\otimes_{A}^2.$$ More concretely, $\mathrm{Sym}_{A}^2 (P) \simeq k\otimes_{k[\Sigma_2] }(QP\otimes_{\mathrm{N}(A)} QP)$.\\ 
Note that, again because $2$ is invertible in $k$, $\Sigma_2$-coinvariants and $\Sigma_2$-invariants are isomorphic, and we might as well have taken invariants in the above construction.\\ Notice that equivalently, we have that $\mathrm{Sym}_{A}^2 (P)$ is defined via the following (usual !) pushout in $A-\dgmod$: $$\xymatrix{QP\otimes_{\mathrm{N}(A)} QP \ar[d]_-{\textrm{1}-\sigma} \ar[r] & 0 \ar[d] \\ QP\otimes_{\mathrm{N}(A)} QP \ar[r] & \mathrm{Sym}_{A}^2 (P) }$$ where $\sigma$ is the symmetry map $p_i \otimes p_j \mapsto (-1)^{ij} p_j \otimes p_i$, for $p_i \in P^i$, and $p_j \in P^j$. \\ An analogous construction, with the properly modified $\Sigma_2$-action, defines the $\infty$-functor $\bigwedge^2_A: A-\underline{\dgmod} \longrightarrow A-\underline{\dgmod}$.



\end{itemize}

\begin{df}\label{dq} Let $C \in A-\dgmod$, and $n\in \mathbb{Z}$.
\begin{itemize}
\item The \emph{space of $n$-shifted derived quadratic forms on} $C$ is the mapping space $$\mathbf{QF}_{A}(C;n):= \mathrm{Map}_{A-\dgmod}(\mathrm{Sym}_{A}^2(C), A[n]) \simeq \mathrm{Map}_{A-\idgmod}(\mathrm{Sym}_{A}^2(C), A[n]) \,\,;$$
\item The \emph{set of $n$-shifted derived quadratic forms on} $C$ is $$\mathrm{QF}_{A}(C;n):= \pi_{0} \mathbf{QF}_{A}(C;n) $$ of connected components of $\mathbf{QF}_{A}(C;n)$, and an \emph{$n$-shifted quadratic form on} $C$ is by definition an element $q\in \mathrm{QF}_{A}(C;n)$.
\item The \emph{space} $\mathbf{QF}_{A}^{\textrm{nd}}(C;n)$ \emph{of $n$-shifted derived non-degenerate quadratic forms on} $C$ is defined by the following homotopy pullback diagram of simplicial sets
$$\xymatrix{\mathbf{QF}_{A}^{\textrm{nd}}(C;n) \ar[d] \ar[r] & \mathbf{QF}_{A}(C;n) \ar[d] \\ [\,\mathrm{Sym}_{A}^2(C), A[n]\,]^{\textrm{nd}} \ar[r] & [\, \mathrm{Sym}_{A}^2(C), A[n]\,] }$$ where $[-,-]$ denotes the hom-sets in the homotopy category of $A-\dgmod$, and $ [\,\mathrm{Sym}_{A}^2(C), A[n]\,]^{\textrm{nd}}$ is the subset of $ [\,\mathrm{Sym}_{A}^2(C), A[n]\,]$ consisting of maps $v: \mathrm{Sym}_{A}^2(C) \to A[n]$ such that the associated, adjoint map $v^{\flat}: C \to C^{\vee}[n]$ is an isomorphism in $\ho (A-\dgmod)$ (where $C^{\vee}:= \mathbb{R}\underline{\mathrm{Hom}_{A}}(C,A)$ is the derived dual of $C$ over $A$).
\item The \emph{set} $\mathrm{QF}_{A}^{\textrm{nd}}(C;n)$ \emph{of $n$-shifted derived non-degenerate quadratic forms on} $C$ is the set $$\mathrm{QF}_{A}^{\textrm{nd}}(C;n):= \pi_{0} \mathbf{QF}_{A}(C;n)^{\textrm{nd}} $$ of connected components of $\mathbf{QF}(C;n)_{A}^{\textrm{nd}}$, and an \emph{$n$-shifted non-degenerate quadratic form on} $C$ is by definition an element $q\in \mathrm{QF}_{A}(C;n)^{\textrm{nd}}$.
\item an $n$\emph{-shifted derived} (resp. \emph{non-degenerate}) \emph{quadratic complex} is a pair $(C,q)$ where $C\in A-\dgmod$ and $q\in \mathrm{QF}_{A}(C;n)$ (resp. $q\in \mathrm{QF}_{A}^{\textrm{nd}}(C;n)$)

\end{itemize} 
\end{df}

Note that, by definition, an $n$-shifted quadratic form $q$ on $C$ is a map $q: \mathrm{Sym}_{A}^2(C) \to A[n]$ in the homotopy category $\ho(A-\dgmod)$.\\
When working over a fixed base $A$, we will often write $\mathbf{QF}(C;n)$ for $\mathbf{QF}_{A}(C;n)$, $\mathbf{QF}(C;n)^{\textrm{nd}}$ for $\mathbf{QF}_{A}(C;n)^{\textrm{nd}}$, and so on.

\begin{rmk} \emph{Note that Definition \ref{dq} makes perfect sense if the cdga $\mathrm{N}(A)$ is replaced by an arbitrary unbounded commutative differential graded algebra (over $k$). Note also that, if $A=k$, the constant simplicial algebra with value $k$, then $\mathrm{N}(A)= k$, and $A-\dgmod = \co$. However, the only nicely behaved situation is when $C$ is a connective dg module (i.e. with vanishing cohomologies in strictly positive degrees) over a non-positively graded cdga. }
\end{rmk}


\begin{exs}\label{example}\emph{
\begin{enumerate}
\item For any $n\in \mathbb{Z}$, any $C \in A-\dgmod$ is an $n$-shifted derived quadratic complex when endowed with the zero derived quadratic form, and for the loop space $\Omega_{0}(\mathbf{QF}(C;n))$ at $0$ of $\mathbf{QF}(C;n))$, we have  an isomorphism in $\ho (\ssets)$ $$\Omega_{0}(\mathbf{QF}(C;n)) \simeq \mathbf{QF}(C;n-1),$$ since $\Omega_{0}\mathrm{Map}_{A-\dgmod}(E, F) \simeq \mathrm{Map}_{A-\dgmod}(E, F[-1]) $, for any $E, F \in A-\dgmod$ 
\item Suppose that $C$ is a connective $m$-connected cochain complex over $k$, for $m>0$ (i.e. $\mathrm{H}^i (C)=0$ for $i>0$ and $-m \leq i\leq 0$). Then, if $n> m+2$, there are no non-zero derived $n$-shifted quadratic forms on $C$ over $k$. In other words, $\mathbf{QF}_{k}(C;n)$ is connected, for any $n>m+2$. This follows immediately from the fact that, under the connectivity hypotheses on $C$, $\mathrm{Sym}_{k}^2(C)$ is $(m+2)$-connected. 
\item If $(C,q,n=0)$ is a derived non-degenerate quadratic complex in $\com$ over $A=k$, then $C$ is discrete (i.e. cohomologically concentrated in degree $0$).
\item If $n,m \in \mathbb{Z}$, and $C \in A-\dgmod $ is a perfect complex, the associated \emph{$(n,m)$-hyperbolic space of} $C$  is the derived $(n+m)$-shifted non-degenerate quadratic complex $\mathsf{hyp}(C;n,m):=(C[n]\oplus C^{\vee}[m], q^{\textrm{hyp}}_{n,m})$, where $q^{\textrm{hyp}}_{(n,m)}$ is given by the composite $$\xymatrix{\mathrm{Sym}^2(C[n]\oplus C^{\vee}[m]) \ar[r]^-{\textrm{pr}} & C[n]\otimes C^{\vee}[m] \simeq (C\otimes C^{\vee})[m+n] \ar[rr]^-{\textrm{ev}[n+m]} & & A[m+n]}$$ where $pr$ denotes the canonical projection and $ev$ the canonical pairing.
Note that, if $A=k$, here we may also take $C=V[0]$, with $V$ a projective finitely generated $k$-module. 
\item Let $A=k$ and $(V,\omega)$ be a symplectic projective and finitely generated $k$-module. Since $\mathrm{Sym}^2 (V[\pm 1]) \simeq \wedge^2 V [\pm 2]$, we get induced derived quadratic non-degenerate structures on $V[1]$ (with shift $-2$), and on $V[-1]$ (with shift $2$). More generally, for any $m \in \mathbb{Z}$, $\mathrm{Sym}^2 (V[2m+1])$ have  an induced $(-4m-2)$-shifted quadratic non-degenerate structure.
\item Let $M=M^n$ be a compact oriented manifold of dimension $n$, $k=\mathbb{F}$ a field of characteristic zero, and let $C:=C^{\bullet}(M;\mathbb{F})$ be cofibrant cdga model, with multiplication $\mu$, for its $\mathbb{F}$-valued singular cochain complex (see, e.g. \cite{su}). Consider the composite $$\xymatrix{q: C\otimes_{\mathbb{F}} C \ar[r]^-{\mu} & C \ar[rr]^-{-\cap [M]} & & \mathbb{F}[n]}.$$ Then, $q$ is obviously symmetric, and it is also non-degenerate (in the derived sense) by Poincar\'e duality. Hence $q$ is an $n$-shifted non-degenerate derived quadratic form on $C=C^{\bullet}(M;\mathbb{F})$, i.e. an element in $\mathrm{QF}^{\mathrm{nd}}_{\mathbb{F}}(C;n)$. 
\item Other symmetric dualities in topology give rise to derived quadratic forms (see, e.g. \cite[Thm 2.5.2, Thm. 3.1.1]{malm}).
\end{enumerate}}
\end{exs}

Let $n \in \mathbb{Z}$. If $\varphi :A \to B'$ a morphism in $\salg$, and $(C,q)$ is a derived $n$-shifted quadratic complex over $A$, the derived base change complex $\varphi^*{C}=C\otimes_{A}B \in \ho (B-\mathsf{dgmod})$ comes naturally endowed with a derived $n$-shifted quadratic form $$\xymatrix{\varphi^{*}q: \mathrm{Sym}_{B}^{2}(C\otimes_{A}B) \ar[r]^-{\sim} & \mathrm{Sym}_{A}^{2}(C)\otimes_{A}B \ar[r]^{q\otimes\mathrm{id}} & A[n]\otimes_{A}B\simeq B[n]}.$$

\begin{df}\label{basechange} Let $n\in \mathbb{Z}$, $(C,q)$ a derived $n$-shifted quadratic complex over $A$, and $\varphi :A \to B$ a morphism in $\salg$. The derived $n$-shifted quadratic complex $(\varphi^{*}C, \varphi^* q)$ over $B$ is called the \emph{base-change} of $(C,q)$  along $\varphi$.
\end{df}

More generally, a morphism $A\to B$ in $\salg$ defines base change maps in $\ho (\ssets)$ $$\varphi^*: \mathbf{QF}_{A}(C;n)\longrightarrow \mathbf{QF}_{B}(\varphi^{*}C;n),$$ and 
$$\varphi^*: \mathbf{QF}_{A}^{\textrm{nd}}(C;n)\longrightarrow \mathbf{QF}_{B}^{\textrm{nd}}(\varphi^{*}C;n).$$

Let $(C_1,q_1)$ and $(C_2,q_2)$ be two $n$-shifted derived quadratic complexes over $A$. Since $$\mathrm{Sym}_{A}^2(C_1\oplus C_2) \simeq \mathrm{Sym}_{A}^2(C_1) \oplus \mathrm{Sym}_{A}^2 (C_2) \oplus (C_1 \otimes_{A}C_2) ,$$ we have a canonical projection $$\pi_{\oplus}: \mathrm{Sym}_{A}^2(C_1\oplus C_2) \longrightarrow \mathrm{Sym}_{A}^2(C_1) \oplus \mathrm{Sym}_{A}^2 (C_2)$$ in $\ho (A-\dgmod)$, that can be used to give the following

\begin{df}\label{sum} Let $(C_1,q_1)$ and $(C_2,q_2)$ be two $n$-shifted derived quadratic complexes over $A$. The \emph{orthogonal sum} $(C_1\oplus C_2, q_1 \perp q_2)$ is the derived $n$-shifted quadratic complex over $A$ where $q_1 \perp q_2$ is defined by the composition $$\xymatrix{\mathrm{Sym}_{A}^2(C_1\oplus C_2) \ar[r]^-{\pi_{\oplus}} & \mathrm{Sym}_{A}^2(C_1) \oplus \mathrm{Sym}_{A}^2 (C_2) \ar[r]^-{q_1\oplus q_2} & A[n]\oplus A[n] \ar[rr]^-{\textrm{sum}} & & A[n]}.$$
\end{df}

\begin{df}\label{res} If $C_1 \in A-\dgmod$, $(C_2,q_2)$ is a derived $n$-shifted quadratic complexes over $A$, and $f:C_1 \to C_2$ is a map in $\ho (A-\dgmod)$, then the composite $$\xymatrix{\mathrm{Sym}_{A}^{2}C_1 \ar[rr]^-{\mathrm{Sym}^{2} f} & & \mathrm{Sym}_{A}^{2}C_2 \ar[r]^-{q_2} & A[n]}$$ defines an $n$-shifted quadratic form on $C_1$, that we denote by $f^*q_2$. $f^*q_2$ is called the \emph{pull-back} or \emph{restriction} of $q_2$ along $f$.
\end{df}

In the present derived setting, the concept of $f$ being an isometry, is not a property of but rather a datum on $f$. More precisely, we give the following

\begin{df}\label{d3} \begin{itemize}
\item The \emph{space of derived isometric structures} on a map $f:(C_1,q_1)\to (C_2,q_2)$, between two $n$-shifted derived quadratic complexes over $A$, is by definition the space $$\mathsf{Isom}(f;(C_1,q_1),(C_2,q_2)):= \mathrm{Path}_{q_1, f^*q_2}(\mathbf{QF}_{A}(C_1;n)).$$ 
\item The space $\mathsf{Isom}(f;(C_1,0),(C_2,q_2))$ is called the \emph{space of derived null-structures} on $f$.
\item A \emph{derived isometric structure} on a map $f:(C_1,q_1)\to (C_2,q_2)$, between two $n$-shifted derived quadratic complexes, is an element in $\pi_{0}(\mathsf{Isom}(f;(C_1,q_1),(C_2,q_2)))$. 
\item A \emph{derived null structure} on a map $f:(C_1,0)\to (C_2,q_2)$, between two $n$-shifted derived quadratic complexes, is an element in $\pi_{0}(\mathsf{Isom}(f;(C_1,0),(C_2,q_2)))$. 
\end{itemize}
\end{df}

\begin{rmk}\textsf{Derived symplectomorphism structures.} \emph{ A similar idea as in the definition above yields a natural notion of \emph{derived symplectomorphism structure} in the theory developed in \cite{ptvv}. If $ \omega $ and $\omega'$  are derived $n$-shifted symplectic forms on a derived stack $X$, the space of \emph{derived symplectic equivalences} between $\omega$ and $\omega'$ is the space $$\mathsf{SymplEq}(X; \omega, \omega' ; n):= \mathrm{Path}_{\omega, \omega'}(\textbf{Sympl}(X;n)).$$ A \emph{derived symplectic equivalence} between $\omega$ and $\omega'$ is then an element $\gamma_{\omega,\omega'} \in \pi_0 \mathsf{SymplEq}(X; \omega, \omega' ; n)$.\\
Let now $f:X_1\to X_2$ be a map, and $\omega_i \in \textbf{Sympl}(X_i;n)$, $i=1,2$. Then, in general, $f^*\omega_2 \in \mathcal{A}^{2,\textrm{cl}}(X_1;n)$ (\cite[\S 2.2]{ptvv}); if moreover $f^* \omega_2 \in \textbf{Sympl}(X_1;n)$, then the  space of \emph{derived symplectomorphism structures} on $f$ is the space $$\mathsf{SymplMor} (f; (X,\omega_1), (X_2,\omega_2)):= \mathsf{SymplEq}(X; \omega_1, f^* \omega_2 ; n),$$ and a \emph{derived symplectomorphism structure} on $f$ is a derived symplectic equivalence $$\gamma_{\omega_1,f^* \omega_2} \in \pi_0 \mathsf{SymplEq}(X_1 ; \omega_1, f^* \omega_2 ; n).$$}
\end{rmk}

\bigskip 

We can now define a notion of \emph{lagrangian structure} (i.e. non-degenerate null structure) on a map of $A$-dg-modules $f: P \to (C,q)$ where $(C,q)$ is a quadratic $A$-dg-module. Let $K$ be the co-cone of $f$, then the composite $$K\rightarrow P \rightarrow C$$ comes with a natural homotopy to the zero map $0: K \to C$, and thus the composite $$\xymatrix{ q_{K}:  K \otimes P \ar[r] & P\otimes P \ar[r]^-{f\otimes f}  & C\otimes C \ar[r]^-{q} & A[n]}$$ acquires an induced homotopy $h_K$ to the zero map. If now $\gamma $ is a derived null structure $\gamma$ on $f: P \to (C,q)$, $\gamma$ induces, by composition, a homotopy $h_{\gamma}$ between $q_{K}$ and the zero map $K\otimes P \to A[n]$. By composing the homotopies (paths) $h_K$ and $h_{\gamma}$, we obtain a loop $\theta_{\gamma}$ at $0$ in the space $\mathrm{Map}_{A-\dgmod}(K\otimes P, A[n])$, i.e. an element, still denoted by $\theta_{\gamma}$, in  $\pi_0 (\mathrm{Map}_{A-\dgmod}(K\otimes P, A[n-1]))\simeq \mathrm{Hom}_{\mathrm{D}(A)}(K\otimes P, A[n-1])$.

\begin{df}\label{d4} \begin{itemize}
\item A derived null structure $\gamma$ on $f: P \to (C,q)$ is said to be \emph{non-degenerate} or \emph{lagrangian} if the induced map $\theta_{\gamma}^{\flat}: K \to P^{\vee}[n-1]$ is an isomorphism in the derived category $\mathrm{D}(A)$.
\item The space $\mathsf{Lag}(f: P\to (C,q))$ of lagrangian structures on the map $f: P\to (C,q)$ is the subspace of $\mathsf{Isom}(f;(P,0), (C,q))$ union of connected components corresponding to lagrangian null-structures.
\end{itemize}
\end{df}

\begin{rmk}
\emph{For $A=k$, $n=0$, $P$ and $C$ both concentrated in degree $0$, and $f:P \hookrightarrow C$ a monomorphism, it is easy to see that the space of lagrangian structures on $f$ is either empty or contractible; it is contractible iff $q$ restricts to zero on $P$, and $C/P \simeq P^{\vee}$, i.e. we recover the usual definition of lagrangian subspace $P$ of a quadratic $k$-space $(C,q)$.}
\end{rmk}

\begin{rmk}
\emph{\textsf{Lagrangian correspondences.} Lagrangian structures can be used in order to define an $\infty$-category $A-\underline{\mathsf{QMod}}_{n}^{\mathsf{Lag}}$ of derived $n$-shifted quadratic modules with morphisms given by  \emph{lagrangian correspondences}. Roughly speaking, the objects are quadratic $n$-shifted $A$-modules $(C,q)$, and the space of maps $(C_1,q_1) \to (C_2, q_2)$ is given by the disjoint union $\coprod_{f}  \mathsf{Lag}(f: P\to (C_1\oplus C_2) (q_1 \perp (-q_2))$ (with composition given by composition of correspondences). This $A-\underline{\mathsf{QMod}}_{n}^{\mathsf{Lag}}$ can be shown to be a symmetric monoidal $\infty$-category under direct sum (see also \cite[\S 7.4]{tt}), and leads to a derived ($n$-shifted) Witt symmetric monoidal $1$-category $\mathrm{h}(A-\underline{\mathsf{QMod}}_{n}^{\mathsf{Lag}})$ of $A$, and further down, to a derived ($n$-shifted) Witt group $\pi_0(\mathrm{h}(A-\underline{\mathsf{QMod}}_{n}^{\mathsf{Lag}}))$ of $A$, with sum given by the orthogonal sum (for $k=A$ a field, $n=0$, and considering only complexes concentrated in degree $0$, we get back the classical Witt group $\mathrm{W}(k)$).  We will not use this sketchy construction further in this paper.   }
\end{rmk}


\section{Derived Clifford algebra of a derived quadratic complex} 
For \emph{even} shifted derived quadratic complexes, it is possible to define a derived version of the Clifford algebra.\\
\subsection{Derived Clifford algebra of a derived quadratic complex} \label{dercliff}

Let $n\in \mathbb{Z}$, and $(C,q)$ be derived $2n$-shifted quadratic complex over $A \in \salg$. We denote by $2\widetilde{q}: C\otimes^{}_{A}C \to A[2n]$ the composite of $2q$ with the canonical map $C\otimes_{A}C  \to \mathrm{Sym}_{A}^2(C)$.\\
If now $B\in A-\dga$, $\widetilde{q}$ induces a map $\widetilde{q}_{B}: * \to \mathrm{Map}_{A-\dgmod}(C\otimes^{}_{A}C, B[2n])$, and the rule $$(\varphi:C\to B[n]) \longmapsto \\ \xymatrix{(C\otimes^{}_{A}C \ar[rr]^-{(\varphi, \varphi \circ \sigma)} & &(B[n]\otimes_{A}^{}B[n]) \oplus (B[n]\otimes^{}_{A}A[n])}$$ $$\xymatrix{ & \ar[r]^-{\mu \oplus \mu} & B[2n]\oplus B[2n] \ar[r]^-{+} & B[2n])} $$ determines a map $$s_B: \mathrm{Map}_{A-\dgmod}(C,B[n]) \to \mathrm{Map}_{A-\dgmod}(C\otimes^{}_{A}C,B[2n]).$$
Here, $$\sigma: C\otimes^{}_{A}C \to C\otimes^{}_{A}C : x\otimes y \mapsto (-1)^{|x||y|} y\otimes x$$ is the Koszul sign involution, $\mu$ denotes the multiplication map on $B$, and $+$ the sum in $B$, so, essentially, the image of $\varphi$
sends $x \otimes y$ to $\varphi(x)\varphi(y) + (-1)^{|x||y|} \varphi(y)\varphi(x)$.\\
By using these maps we may define the derived Clifford algebra functor associated to the derived $2n$-shifted quadratic space $(C,q)$, as $$\underline{\mathbf{Cliff}}(C,q,2n): A-\dga \longrightarrow \ssets \,\, : \, B \longmapsto \underline{\mathbf{Cliff}}(C,q,2n)(B)$$ where $\underline{\mathbf{Cliff}}(C,q,2n)(B)$ is defined by the following homotopy pull back in $\ssets$ $$\xymatrix{\underline{\mathbf{Cliff}}(C,q,2n)(B) \ar[r] \ar[d] & \mathrm{Map}_{A-\dgmod}(C,B[n]) \ar[d]^-{s_B} \\ \textrm{*} \ar[r]_-{\widetilde{q}_{B}} &  \mathrm{Map}_{A-\dgmod}(C\otimes^{}_{A}C,B[2n])}$$

\begin{prop}\label{Cliff} The functor $\underline{\mathbf{Cliff}}(C,q,2n)$ is homotopy co-representable, i.e. there exists a well defined $\mathbf{Cliff}_{A}(C,q,2n) \in \ho (A-\dga)$ and a canonical isomorphism in $\ho (\ssets)$ $$\underline{\mathbf{Cliff}}(C,q,2n)(B)\simeq \mathrm{Map}_{A-\dga}(\mathbf{Cliff}_{A}(C,q,2n),B).$$
\end{prop}

\noindent \textbf{Proof.}   Since the notion of ideal is not well-behaved in derived geometry, we need to reformulate the existence in a homotopical meaningful way. This leads us to the following construction. Let $\mathsf{Free}_{A}: A-\dgmod \to A-\dga$ be the left derived functor of  the free dga-functor (Quillen left adjoint to the forgetful functor). 
Then $\mathbf{Cliff}(C,q;2n)$ is defined by the following homotopy push-out square in $A-\dga$ $$\xymatrix{\mathsf{Free}_{A}(C\otimes_{A}C[-2n]) \ar[r]^-{u} \ar[d]_-{t_q} & \mathsf{Free}_{A}(C[-n]) \ar[d] \\ A \ar[r] & \mathbf{Cliff}(C,q,2n) }$$ where 
\begin{itemize}
\item $t_q$ is induced, by adjunction, by the map $2\widetilde{q}[-2n]: C[-n]\otimes^{}_{A}C[-n] \to A[2n-2n]=A$
\item $u$ is induced, by adjunction, by the map $$\xymatrix{C[-n]\otimes^{}_{A}C[-n] \ar[r]^-{(id, \sigma)}  &(C[-n]\otimes_{A}^{}C[-n]) \oplus (C[-n]\otimes^{}_{A}C[-n])  \ar[r]^-{+}  & }$$ $$\xymatrix{  \ar[r]^-{+} & C[-n]\otimes_{A}^{}C[-n] \ar[r] & \mathsf{Free}_{A}(C[-n]).}$$
\end{itemize}
\hfill $\Box$ \\
 
 \begin{df} The dga $\mathbf{Cliff}_{A}(C,q;2n)$, defined up to isomorphism in $\ho (A-\dga)$, is called the \emph{derived Clifford algebra} of the derived $2n$-shifted quadratic space $(C,q)$.
 \end{df}
 
When the base simplicial algebra $A$ is clear from the context, we will simply write $\mathbf{Cliff}(C,q,2n)$ for $\mathbf{Cliff}_{A}(C,q,2n)$.

\begin{rmk}\label{classicalCL} \emph{To understand classically the proof of Proposition \ref{cliff}, observe that the classical Clifford algebra of a quadratic $k$-module $(V,Q)$ is defined as the quotient $$\frac{\mathrm{T}_{k}(V)}{\mathrm{I}:=<x\otimes y + y\otimes x -2Q(x,y)>}\,\, ,$$ $\mathrm{T}_{k}(-)$ denoting the tensor $k$-algebra functor. Then it is easy to verify that the following square is a (strict) push-out in the category of $k$-algebras $$\xymatrix{\mathrm{T}_{k}(V\otimes_k V) \ar[r]^-{u} \ar[d]_-{t} & \mathrm{T}_{k}(V) \ar[d] \\ k \ar[r] & \mathrm{T}_{k}(V)/\mathrm{I} }$$ where $u$ is defined by $u(x):=x\otimes y + y\otimes x$, and $t$ by $t(x\otimes y):= 2Q(x,y)$. For a more detailed study of the relation between the classical and the derived Clifford algebra when they both apply, see Section \ref{comparison}.}
\end{rmk}

\begin{rmk}\label{finitetype} \emph{Note that it follows immediately from the description of $\mathbf{Cliff}(C,q,2n)$ in the proof of Proposition \ref{cliff}, that $\mathbf{Cliff}(C,q,2n)$ is a dg algebra \emph{homotopically of finite type} over $A$ whenever $C$ is a perfect $A$-dg module. }
\end{rmk}

\begin{prop}\label{pres}  Let $n \in \mathbb{Z}$, $f: C_1 \to C_2$ be a map in $A-\dgmod$, and $q_2$ a derived $2n$-shifted quadratic form on $C_2$ over $A$. Then there is a canonical map in $\ho (\dga)$ $$\mathbf{Cliff}_{A}(C_1,f^*q_2,2n) \longmapsto \mathbf{Cliff}_{A}(C_2,q_2,2n)$$ where $f^*q_2$ is the pull-back quadratic form of Definition \ref{res}.

\end{prop}

\noindent \textbf{Proof.}   We use the proof of Proposition \ref{cliff}. By definition of $f^*q_2$, there is a map (in the homotopy category of diagrams of that shape in $\dga$) from the diagram $$\xymatrix{\mathsf{Free}_{A}(C_1\otimes_{A}C_1[-2n]) \ar[r]^-{u} \ar[d]_-{t_{f^* q_2}} & \mathsf{Free}_{A}(C_1[-n])  \\ A &  }$$ to the diagram $$\xymatrix{\mathsf{Free}_{A}(C_2\otimes_{A}C_2[-2n]) \ar[r]^-{u} \ar[d]_-{t_{q_2}} & \mathsf{Free}_{A}(C_2[-n]) \\ A &  }.$$ By definition of 
homotopy push-out, we get the induced map.  \hfill $\Box$ \\
 




Note that, by composition with the natural adjunction map of complexes $C[-n] \rightarrow \mathsf{Free}_{A}(C[-n])$, the derived Clifford algebra $\mathbf{Cliff}_{A}(C,q,2n)$ of  the derived $2n$-shifted quadratic complex $(C,q)$, comes equipped with a natural map of $A$-dg modules $$C[-n] \longrightarrow \mathbf{Cliff}(C,q,2n).$$
Using this, and the universal property of the derived Clifford algebra, we get the following

\begin{prop}\label{p1} Let $n \in \mathbb{Z}$. If $(C_1,q_1)$ and $(C_2,q_2)$ are $2n$-shifted derived quadratic complexes over $A$, $f:C_1 \to C_2$ is a map in $\ho (A-\dgmod)$ and $\gamma \in  \pi_{0}(\mathsf{Isom}(f;(C_1,q_1),(C_2,q_2)))$ is a derived isometric structure on $f$, then there is a canonical induced map in $\ho (A-\dga)$ $$f_{\gamma}: \mathbf{Cliff}(C_1,q_1,2n) \longrightarrow \mathbf{Cliff}(C_2,q_2,2n).$$ If moreover $f$ is a quasi-isomorphism, $f_{\gamma}$ is an isomorphism.
\end{prop}

\noindent \textbf{Proof.}  Let us consider the following two maps $$\xymatrix{h:  \mathsf{Free}_{A}(C_1\otimes^{}_{A}C_1[-2n]) \ar[r]^-{u_{C_1}} & \mathsf{Free}_{A}(C_1[-n]) \ar[rr]^-{\mathsf{Free}(f[-n])} &  & \mathsf{Free}_{A}(C_2[-n]) \ar[r] & \mathbf{Cliff}(C_2,q_2,2n),  }$$ and $$\xymatrix{g_1:  \mathsf{Free}_{A}(C_1\otimes^{}_{A}C_1[-2n]) \ar[r]^-{\tilde{q_1}} & A \ar[r] & \mathbf{Cliff}(C_2,q_2,2n),  }$$

 By definition of $\mathbf{Cliff}(C_1,q_1,2n)$ as a homotopy push-out (proof of Proposition \ref{cliff}),  it is enough to show that our data give a path between $h$ and $g_1$.  To start with, the homotopy push-out defining  $\mathbf{Cliff}(C_2,q_2,2n)$, provides us with a path between the composite maps $$\xymatrix{ & \mathsf{Free}_{A}(C_2\otimes^{}_{A}C_2[-2n]) \ar[r]^-{u_{C_2}} & \mathsf{Free}_{A}(C_1[-n]) \ar[r] & \mathbf{Cliff}(C_2,q_2,2n),  }$$ and 
 $$\xymatrix{& \mathsf{Free}_{A}(C_2\otimes^{}_{A}C_2[-2n]) \ar[r]^-{\tilde{q_2}} & A \ar[r] & \mathbf{Cliff}(C_2,q_2,2n),  }$$
 and, by precomposing with $\mathsf{Free}_{A}(f[-n])$, we get a path $\delta$ between the maps
 $$\xymatrix{h':  \mathsf{Free}_{A}(C_1\otimes^{}_{A}C_1[-2n]) \ar[rr]^-{\mathsf{Free}(f[-n])} & & \mathsf{Free}_{A}(C_2\otimes^{}_{A}C_2[-2n]) \ar[r]^-{u_{C_2}} & \mathsf{Free}_{A}(C_1[-n]) \ar[r] & \mathbf{Cliff}(C_2,q_2,2n),  }$$ and $$\xymatrix{g'_2 : \mathsf{Free}_{A}(C_1\otimes^{}_{A}C_1[-2n]) \ar[rr]^-{\mathsf{Free}(f[-n])} & & \mathsf{Free}_{A}(C_2\otimes^{}_{A}C_2[-2n]) \ar[r]^-{\tilde{q_2}} & A \ar[r] & \mathbf{Cliff}(C_2,q_2,2n). }$$ Now observe that $h' = h$, while $g'_{2}$ is equal to the composite $$\xymatrix{g_2:  \mathsf{Free}_{A}(C_1\otimes^{}_{A}C_1[-2n]) \ar[r]^-{\widetilde{f^* q_2}} & A \ar[r] & \mathbf{Cliff}(C_2,q_2,2n),  }$$ by definition of $\widetilde{f^* q_2}$. Hence, $\delta$ gives us a path between $h$ and $g_2$. We conclude by using the further path between $\tilde{q_1}$ and $\widetilde{f^* q_2}$ (hence between $g_2$ and $g_1$), induced by the derived isometric structure on $f$, i.e. by the path $\gamma$ between $q_{1}$ and $f^* q_2$.

 \hfill $\Box$\\

In the following Proposition, recall the definition of base change for a derived quadratic complex along a base ring morphism (Definition \ref{basechange}).

\begin{prop}\label{p2} Let $n \in \mathbb{Z}$, $(C,q)$ be a $2n$-shifted derived quadratic complex over $A$, and $\varphi :A \to B$ a morphism in $\salg$. Then there is a canonical isomorphism in $\ho (B-\mathsf{dga})$ $$ \mathbf{Cliff}_{A}(C,q,2n)\otimes_{A}B \simeq \mathbf{Cliff}_{B}(\varphi^{*}C,\varphi^{*}q,2n).$$ 
\end{prop}

\noindent \textbf{Proof.}  This follows immediately from the proof of Proposition \ref{cliff}, and the observation that there exists a natural isomorphism $\mathsf{Free}_{A}(E)\otimes_{A} B \simeq \mathsf{Free}_{B}(E\otimes_{A} B)$ in $\ho (B-\mathsf{dga})$, for any $E \in A-\dgmod$.

 \hfill $\Box$ \\

As shown in Appendix \ref{app}, if $n \in \mathbb{Z}$, and $(C,q)$ be a $2n$-shifted derived quadratic complex over $A$, then $\mathbf{Cliff}_{A}(C,q,2n)$ admits a natural $\mathbb{Z}/2$-weight grading, i.e it is naturally an object $\mathbf{Cliff}^{w}_{A}(C,q,2n)\in \ho(A-\dga^{w})$ (see the Appendix for the notations). Then, a slight modification of classical proof of the corresponding classical statement (see e.g. \cite[Prop. 2.2.1]{mire}), yields the following

\begin{prop}\label{graded} Let $n \in \mathbb{Z}$, $(C_i,q_i)$ be $2n$-shifted derived quadratic complexes over $A$, and $(C_1 \oplus C_2, q_1 \perp q_2)$ be the corresponding orthogonal sum (Def. \ref{sum}). Then there is a canonical isomorphism in $\ho (A-\mathsf{dga}^{w})$ $$ \mathbf{Cliff}^{w}_{A}(C_1 \oplus C_2,q_1 \perp q_2, 2n)\simeq \mathbf{Cliff}^{w}_{A}(C_1 ,q_1, 2n) \otimes^{w}_{A}\mathbf{Cliff}^{w}_{A}(C_2, q_2, 2n),$$ where $\otimes^{w}_{A}$ denotes the derived tensor product of $\mathbb{Z}/2$-dg algebras over $A$ (see \S \ref{app}). 
\end{prop}

\begin{rmk}\emph{\textsf{From derived Clifford algebras to derived Azumaya algebras ?} $\,$ 
Classically, Clifford algebras are Azumaya algebras: this is literally true when the rank of the quadratic module is even, true for its even graded piece in the odd rank case, and always true if we replace Azumaya algebras by
$\mathbb{Z}/2$-graded Azumaya algebras (i.e. those classified by the Brauer-Wall group).  Recently, To\"en introduced the notion of derived Azumaya algebras over a (derived) stack, and proved e.g. that they are classified, up to Morita equivalence, by $H_{\textrm{ét}}^1 (X,\mathbb{Z})\times H_{\textrm{ét}}^2 (X,\mathbb{G}_{m} )$ on a quasi-separated, quasi-compact scheme $X$ (\cite[Cor. 3.8]{deraz}).\\ 
As suggested by To\"en it might be interesting to investigate the exact relation  between derived Clifford algebras and derived Azumaya algebras. The first problem we meet in such a comparison is the presence of a non-zero shift in the derived quadratic structure. There is no evident place for a shift in the current definition of derived Azumaya algebras, and it is not clear to us how one could modify the definition of a derived Azumaya algebra in order to accommodate such a shift. If we limit ourselves to the $0$-shifted case, one can give at least one (admittedly not very relevant). Let $C$ be a perfect complex of $k$-modules, and consider the associated $0$-hyperbolic space $\mathsf{hyp}(C;n=0,m=0):=(C\oplus C^{\vee}, q^{\textrm{hyp}}_{0,0}; 0)$ of Example \ref{example} (3). Then there is a canonical isomorphism $$\mathbf{Cliff}(\mathsf{hyp}(C;0,0)) \simeq \mathbb{R}\underline{\textrm{Hom}}_{k}(\bigwedge C, \bigwedge C)$$ in $\ho (\dga)$, where $\mathbb{R}\underline{\textrm{Hom}}_{k}$ denotes the derived internal Hom's in $\co$. Thus $\mathbf{Cliff}(\mathsf{hyp}(C;0,0))$ is indeed a derived Azumaya algebra, but its class is obviously trivial in the derived Brauer group introduced by To\"en (\cite[Def. 2.14]{deraz}). Note this case is morally an even rank case. A slightly more interesting statement would be the following (still sticking to the $0$-shifted case). By systematically replacing $\mathbb{Z}/2$-graded dga's (see Appendix \ref{app}) in To\"en's definition of derived Azumaya algebras, we get a notion of $\mathbb{Z}/2$\emph{-graded derived Azumaya algebras}. Note that the extra  $\mathbb{Z}/2$-grading affects both the notion of opposite algebra and of tensor product. As shown in Appendix \ref{app} below, the derived Clifford algebra of derived $2n$-shifted quadratic complex over $k$ is naturally an object in the homotopy category of  $\mathbb{Z}/2$-graded dga's over $k$. So we may formulate the following question that has an affirmative answer in the underived, classical case : is the $\mathbb{Z}/2$-graded  derived Clifford algebra of a $0$-shifted quadratic complex over $k$ a $\mathbb{Z}/2$-graded derived Azumaya algebra?\\ Of course, one could formulate a similar question by replacing $k$ with a base scheme or an arbitrary derived Artin stack. This question is probably not too hard to settle but we do not have neither a proof nor a counterexample at the moment. \\ A related interesting question is to extend To\"en cohomological identification of the derived Brauer group to the corresponding \emph{derived Brauer-Wall group} i.e. the Morita equivalence classes of  derived $\mathbb{Z}/2$-graded Azumaya algebras. }

\end{rmk}

\subsection{Comparison with the classical Clifford algebra}\label{comparison}
In this section, we will work over $A=k$. Given a a classical quadratic (projective and finitely generated) $k$-module $(V,Q)$,  a natural question is how can we get back the usual Clifford algebra $\mathbf{Cliff}^{\textrm{class}}(V,Q )$ from its derived ($0$-shifted) Clifford algebra $\mathbf{Cliff}_{k}(V[0],Q[0],0)$.\\

\begin{prop} Let $(C, q)$ be a derived $0$-shifted quadratic complex in $\com$. Then $(H^0(C), H^0(q))$ is a classical quadratic $k$-module, and there is a natural isomorphism of associative unital algebras $$\mathbf{Cliff}^{\textrm{class}}(H^0(C), H^0(q)) \simeq H^{0}(\mathbf{Cliff}_{k}(C,q,0)).$$
\end{prop}

\noindent \textbf{Proof.}   Let $\mathsf{Free}: \co \to \dga_{k}$ be the free dga-functor (left adjoint to the forgetful functor $\mathsf{For}: \dga_{k} \to \co$), $\mathsf{Free}^{\leq 0}: \com \to \dgam$ the free dga-functor (left adjoint to the forgetful functor $\mathsf{For}^{\leq 0}: \dgam \to \com$), $j_{-}: \com \to \co$ the inclusion functor, $[-]_{\leq 0}: \co \to \com$ its right adjoint (given by the intelligent truncation in degrees $\leq 0$), $i_{-} : \dgam \to \dga_{k}$ the inclusion functor, and $(-)_{\leq 0}: \dga_{k} \to \dgam$ its right adjoint (given by the intelligent truncation in degrees $\leq 0$). Note that the pairs $(j_{-}, [-]_{\leq 0} )$ and $(i_{-}, (-)_{\leq 0})$, are Quillen pairs of type (left, right).\\

\begin{lem} \begin{enumerate}
\item Let $$\xymatrix{A \ar[r] \ar[d]  & B \ar[d]  \\ C \ar[r] & D}$$ be a homotopy push-out diagram in $\dgam$, and $$\xymatrix{i_{-}A \ar[r] \ar[d]  & i_{-}B \ar[d]  \\ i_{-}C \ar[r] & D'}$$ the homotopy push-out in $\dga_{k}$. Then the canonical map $i_{-}D' \to D$ is a weak equivalence.
\item If $E \in \com$, then the canonical map $i_{-}\mathsf{Free}^{\leq 0}(E) \to \mathsf{Free}(j_{-}E)$ is a weak equivalence.
\end{enumerate}
\end{lem}

\noindent \textsl{Proof of Lemma.} Statement 1 follows immediately from the fact that $i_{-}$ is left Quillen, hence preserves homotopy push-outs.\\ To prove 2, we start from the adjunction map $j_{-}E \to \textrm{For}\mathsf{Free}(j_{-}E) $ in $\co$. Since $j_{-}$ is left adjoint, we get a map $E \to [\mathsf{For}\mathsf{Free}(j_{-}E)]_{\leq 0}$ in $\com$. Note that the canonical transformation $[-]_{\leq 0}\circ \mathsf{For} \to \mathsf{For}^{\leq 0} \circ (-)_{\leq 0} $ is an isomorphism, hence we get a map $E \to \mathsf{For}^{\leq 0}((\mathsf{Free}(j_{-}E))_{\leq 0})$ in $\com$. Since $\mathsf{For}^{\leq 0}$ is right-adjoint, we get an induced map $\mathsf{Free}^{\leq 0}(E) \to (\mathsf{Free}(j_{-}E))_{\leq}$ in $\dgam$, and finally the desired canonical map $i_{-}\mathsf{Free}^{\leq 0}(E) \to \mathsf{Free}(j_{-}E)$. Since $E$ and $k$ are connective (i.e. concentrated in non-positive degrees), we have $\mathsf{Free}(j_{-}E) \simeq \oplus_{n\geq 0}E^{\otimes_{k}^{n}}$, and $\otimes_{k}$ is left $t$-exact;  therefore the map $i_{-}\mathsf{Free}^{\leq 0}(E) \to \mathsf{Free}(j_{-}E)$ is indeed a weak equivalence.
 \hfill $\Box$\\

Let us consider now the following homotopy push-out square in $\dgam$ (defining $\mathbf{Cliff}^{\leq 0}(C,q,0)$) $$\xymatrix{\mathsf{Free}^{\leq 0}(C\otimes^{}_{k}C) \ar[r]^-{u^{\leq 0}} \ar[d]_-{t^{\leq 0}} & \mathsf{Free}^{\leq 0}(C) \ar[d] \\ \mathsf{Free}^{\leq 0}(0)=k \ar[r] & \mathbf{Cliff}^{\leq 0}(C,q,0). }$$ By applying the left Quillen functor $i_{-}$, by point 1 in Lemma, we then get a homotopy push-out in $\dga_{k}$

$$\xymatrix{i_{-}\mathsf{Free}^{\leq 0}(C\otimes^{}_{k}C) \ar[r]^-{i_{-}u^{\leq 0}} \ar[d]_-{i_{-}t^{\leq 0}} & i_{-}\mathsf{Free}^{\leq 0}(C) \ar[d] \\ i_{-}\mathsf{Free}^{\leq 0}(0)=k \ar[r] & i_{-}\mathbf{Cliff}^{\leq 0}(C,q,0). }$$

By point 2 in Lemma above, we also have a homotopy push-out diagram  in $\dga_{k}$ $$\xymatrix{\mathsf{Free}(j_{-}C\otimes^{}_{k}j_{-}C) \ar[r]^-{u} \ar[d]_-{t} & i_{-}\mathsf{Free}(j_{-}C) \ar[d] \\ k \ar[r] & i_{-}\mathbf{Cliff}^{\leq 0}(C,q,0), }$$ hence, by definition of  $\mathbf{Cliff}(C,q,0)$, we get an equivalence $ i_{-}\mathbf{Cliff}^{\leq 0}(C,q,0) \simeq  \mathbf{Cliff}(j_{-}C,q,0)$ in $\dga_{k}$ (therefore $\mathbf{Cliff}(j_{-}C,q,0)$ is cohomologically concentrated in non-positive degrees). Since $i_{-}$ preserves weak equivalences, we also get an induced map $$\mathbf{Cliff}(j_{-}C,q,0) \to H^{0}(\mathbf{Cliff}(j_{-}C,q,0))= H^0(\mathbf{Cliff}^{\leq 0}(C,q,0))$$ that is an isomorphism on $H^{0}$. Since $H^0: \dgam \to \mathsf{alg} $ is left Quillen and $H^0 \mathsf{Free}^{\leq 0}(E) \simeq \mathsf{Free}^{0}(H^{0}E) $ (where $E$ is any object in $\com$ and $\mathsf{Free}^{0}$ is the free $k$-algebra functor defined on $k$-modules), by definition of the classical Clifford algebra (see Remark \ref{classicalCL}) we get an isomorphism $$H^0(\mathbf{Cliff}^{\leq 0}(C,q,0)) \simeq \mathbf{Cliff}^{\textrm{class}}(H^0(C), H^0(q)).$$

 \hfill $\Box$ \\
 
\begin{cor} \label{comparazione} If $(V,Q)$ is a quadratic $k$-module, then the canonical map of dga's $$\mathbf{Cliff}_{k}(V[0],Q[0],0) \longrightarrow \mathbf{Cliff}^{\textrm{class}}(V,Q )$$ induces an isomorphism on $H^0$.
\end{cor}

One obvious natural question is now the following: suppose that $(V,Q)$ is a quadratic $k$-module, is the dga $\mathbf{Cliff}_{k}(V[0],Q[0],0)$ (which we know being cohomologically concentrated in non-positive degrees) $0$-truncated (i.e. discrete) ? In other words, we are asking whether the derived Clifford algebra of a classical quadratic $k$-module contains or not strictly more information than its classical Clifford algebra.\\
The following example shows that the answer to this question is that, in general, the derived Clifford algebra of a classical quadratic $k$-module is \emph{not} discrete, so it contains a priori more information than its classical Clifford algebra. The further question, i.e. whether or not this extra information might be relevant to the classical theory of quadratic forms - by giving new invariants or just a reinterpretation of known ones - is interesting but will not be addressed in the present paper.

 \begin{ex}\emph{Let $V=k\oplus k$ be endowed with the quadratic form $Q$ given by the symmetric matrix $(Q_{ij})_{i,j=1,2}$. We will show that $\mathrm{H}^{-1}(\mathbf{Cliff}_{k}(V[0],Q[0],0))\neq 0$ by simply using the symmetry of $(Q_{ij})$. \\Let $$A:= k< x_{ij}\, |\, i,j=1,2> , $$ $$B_1:= k< x_{ij}, y_{ij}\, |\, i,j=1,2> \qquad B:=  k< t_{l} \, |\, l=1,2>$$ where $k<...>$ denotes the free associative graded algebra on the specified generators, and $\deg (x_{ij})= \deg (t_{l})=0$, for any $i,j,l=1,2$, while $\deg (y_{ij})=-1$, for any $i,j=1,2$, and $B_1$ is endowed with the unique differential $d$ such that $d(y_{ij})=x_{ij}$ and making $B_1$ into a differential graded algebra. The map $$t': \, k< x_{ij}\, |\, i,j=1,2>\, \longrightarrow B_{1}\, : \, x_{ij} \longmapsto x_{ij} + 2Q_{ij}$$ is a cofibrant replacement in $k-\dga$ of the map $$t :\, k < x_{ij}\, |\, i,j=1,2>\, \longrightarrow k\, : \, x_{ij}\longmapsto 2Q_{ij}$$ as the factorization of $t$ as $$\xymatrix{ k< x_{ij}\, |\, i,j=1,2>  \ar[r]^-{t'} & B_1 \ar[rr]_{\sim}^-{x_{ij}\mapsto 0} && k }$$ shows. Then $\mathbf{Cliff}_{k}(V[0],Q[0],0)$ can be identified with the strict pushout of the following diagram $$\xymatrix{A \ar[r]^-{u} \ar[d]_-{t'} & B \\ B_1 & }$$ in $k-\dga$ (where $u(x_{ij}):= x_i x_j + x_j x_i$). The pushout of dga's is a rather complicated object to describe concretely, so we use the following trick. \\
Let $C$ be the dga over $k$ whose underlying graded algebra is the free graded algebra on generators $\{y_{ij}\}_{i,j=1,2}$ in degree $-1$, with the unique differential making it into a dga over $k$ such that $d(y_{ij}):= -2Q_{ij}$.
Consider the following maps in $k-\dga$ $$\psi: B=  k< t_{l} \, |\, l=1,2> \longrightarrow C \, : \, x_{l} \longmapsto 0 $$ $$\varphi: B_1= k< x_{ij}, y_{ij}\, |\, i,j=1,2> \longrightarrow C \, : \, x_{ij} \mapsto -2Q_{ij} \,\,, \, y_{ij} \mapsto y_{ij} \, ;$$ it is easy to verify that the two composites $$\xymatrix{A \ar[r]^-{t'} & B_1 \ar[r]^-{\varphi} & C }$$ $$\xymatrix{A \ar[r]^-{u} & B\ar[r]^-{\psi} & C, }$$ coincide. Hence we obtain a canonical map of dga's $$f:\mathbf{Cliff}_{k}(V[0],Q[0],0) \longrightarrow C.$$ Suppose now that $Q_{12}=Q_{21}=0$, and that neither $2Q_{11}$ nor $2Q_{22}$ are invertible in $k$. Then, by a long but straightforward computation, one checks that the element $(y_{ij} - y_{ji})$ is a $(-1)$-cycle and gives a non-zero class in $\overline{\alpha} \in \mathrm{H}^{-1}(C)$, for $i\neq j\,\,$\footnote{There are probably other, non-diagonal, cases for which the statement holds: the general answer boils down to showing that  two matrices of sizes $4\times 16$ and $4 \times 17$ have different ranks over $k$. I hope this explains my choice of a particular solution.}. Note that $\varphi (\alpha:= y_{ij}-y_{ji})=y_{ij}-y_{ji}$. Hence, if  $$p: B_1 \longrightarrow \mathbf{Cliff}_{k}(V[0],Q[0],0) \qquad q:B \longrightarrow \mathbf{Cliff}_{k}(V[0],Q[0],0)$$ are the natural maps, in order to show that  $\mathrm{H}^{-1}(\mathbf{Cliff}_{k}(V[0],Q[0],0))\neq 0$, it will be enough to show that the element $p(\alpha)$ is a $(-1)$-cycle in $\mathbf{Cliff}_{k}(V[0],Q[0],0)$. But this is easy to check, by just using that $Q_{ij}=Q_{ji}$; more precisely: $$d(p(\alpha))= p(x_{ij}-x_{ji})= p((x_{ij}+2 Q_{ij})-(x_{ji} + 2Q_{ji})) = p(t'(x_{ij}-x_{ji}))=$$ 
 $$= q(u(x_{ij}-x_{ji}))= q(x_{i}x_j + x_j x_i -(x_jx_i + x_i x_j))=0.$$ So, the $(-1)$-cycle $p(\alpha)$ in $\mathbf{Cliff}_{k}(V[0],Q[0],0)$, has image a non zero cycle via $$f:\mathbf{Cliff}_{k}(V[0],Q[0],0) \longrightarrow C,$$ hence its class is non-zero in $\mathrm{H}^{-1}(\mathbf{Cliff}_{k}(V[0],Q[0],0))$. Note that this apply also to the case where $k$ is a field with characteristic different from $2$, and $Q=0$. }
\end{ex}

\begin{rmk} \emph{Recently, Bertrand To\"en (private communication), simplified the previous example as follows. Start with $V=k$ a field of characteristic zero, with the zero quadratic form. Then it is easy to verify that $\mathbf{Cliff}_{k}(V[0],0,0)$ is the free graded $k$-algebra generated by $x$ in degree $0$ and $y$ in degree $-1$, with differential defined by $d(y)=x^2$. One then verifies, as in the above example, that this dga has non vanishing $\mathrm{H}^{-1}$, while its truncation is the corresponding exterior algebra i.e. the $k$-algebra of dual numbers $k[x]/x^2$.   }
\end{rmk}

\section{Derived quadratic complexes and derived quadratic stacks}
In this Section, we will follow the theory of derived stacks over $k$ as developed in \cite[\S 2.2]{hagII}.

\subsection{Derived quadratic complexes on a derived stack }

In this section we globalize the notions of Section 1 to derived stacks. For future reference, we give a rather complete treatment, with the possible effect of being a bit pedantic. We apologize to the reader if this is the case. \\

\begin{df}\label{qmodA} Let $A \in \salg$ and $n\in \mathbb{Z}$. \begin{itemize}
\item  Define the $\infty$-category $\underline{\mathsf{QMod}}(A; n)$  via the following fiber product of $\infty$-categories (in the $(\infty,1)$-category of presentable $(\infty,1)$-categories)
$$\xymatrix{\underline{\mathsf{QMod}}(A; n)\ar[rr] \ar[d] & & A-\underline{\dgmod}^{\Delta[1]} \ar[d]^-{\textrm{ev}_0 \times \textrm{ev}_1} \\ A-\underline{\dgmod} \times \Delta[0] \ar[rr]^-{\mathrm{Sym}_{A}^{2} \times A[n]} & & A-\underline{\dgmod} \times A-\underline{\dgmod}  }$$ where $\Delta[i]$ denotes the $i$-simplex (or, equivalently, $\Delta[0]$ is the category with one object and no non-identity maps, and $\Delta[1]$ is the category $\{0 \to 1 \}$).
The $\infty$-category $\underline{\mathsf{QMod}}(A; n)$ is called the \emph{category of of $n$-shifted quadratic complexes over $A$}. \\ The full $\infty$-sub category of non degenerate complexes is denoted by $\underline{\mathsf{QMod}}^{\mathrm{nd}}(A; n)$. 
\item Define the $\infty$-category $\underline{\mathsf{QPerf}}(A; n)$  via the following fiber product of $\infty$-categories (in the $(\infty,1)$-category of presentable $(\infty,1)$-categories)
$$\xymatrix{\underline{\mathsf{QPerf}}(A; n)\ar[rr] \ar[d] & & A-\iperf^{\Delta[1]} \ar[d]^-{\textrm{ev}_0 \times \textrm{ev}_1} \\ A-\iperf \times \Delta[0] \ar[rr]^-{\mathrm{Sym}_{A}^{2} \times A[n]} & & A-\iperf \times A-\iperf  }$$ where $\Delta[i]$ denotes the $i$-simplex (or, equivalently, $\Delta[0]$ is the category with one object and no non-identity maps, and $\Delta[1]$ is the category $\{0 \to 1 \}$).
The $\infty$-category $\underline{\mathsf{QPerf}}(A; n)$ is called the \emph{category of of $n$-shifted quadratic perfect complexes over $A$}. \\ The full $\infty$-sub category of non degenerate (perfect) complexes is denoted by $\underline{\mathsf{QPerf}}^{\mathrm{nd}}(A; n)$. 
\end{itemize}
\end{df}

For any morphism $\varphi: A\to B$ in $\salg$, there are induced base-change $\infty$-functors $$\varphi^* : \underline{\mathsf{QMod}}(A; n) \longrightarrow \underline{\mathsf{QMod}}(B; n),$$
$$\varphi^* : \underline{\mathsf{QPerf}}(A; n) \longrightarrow \underline{\mathsf{QPerf}}(B; n),$$
 so that we get (as explained in \cite[\S 1]{chern}) cofibered $\infty$-categories $\underline{\mathsf{QMod}}(n)$ and $\underline{\mathsf{QPerf}}(n)$ over $\salg$,  whose associated $\infty$-functors will be denoted $$\underline{\mathcal{Q}\mathsf{Mod}} (n): \salg \longrightarrow \infty-\textsf{Cat}, $$ $$\underline{\mathcal{Q}\mathsf{Perf}} (n): \salg \longrightarrow \infty-\textsf{Cat}.$$  Since the fiber product of (derived) stacks is a (derived) stack, we have that both $\underline{\mathcal{Q}\mathsf{Mod}} (n)$ and  $\underline{\mathcal{Q}\mathsf{Perf}} (n)$ are derived stacks, with respect to the derived étale topology, with values in $\infty$-categories, called the \emph{stack of derived $n$-shifted quadratic complexes}, and the \emph{stack of derived $n$-shifted perfect quadratic complexes}, respectively.\\  The \emph{stack of derived $n$-shifted non-degenerate quadratic complexes} is the sub-stack $\underline{\mathcal{Q}\mathsf{Mod}}^{\mathrm{nd}} (n)$ obtained by working with $\underline{\mathsf{QMod}}^{\mathrm{nd}}(A; n)$ instead of $\underline{\mathsf{QMod}}(A; n)$.\\The \emph{stack of derived $n$-shifted non-degenerate perfect quadratic complexes} is the sub-stack $\underline{\mathcal{Q}\mathsf{Perf}}^{\mathrm{nd}} (n)$ obtained by working with $\underline{\mathsf{QPerf}}^{\mathrm{nd}}(A; n)$ instead of $\underline{\mathsf{QPerf}}(A; n)$.\\
Equivalently, the assignments $\lqcoh: \salg \ni A \mapsto A-\idgmod(X)$, and  $\lperf : \salg \ni A \mapsto A-\iperf$ also define derived stacks, with respect to the derived étale topology, with values in $\infty$-categories, and we have pullback diagrams (in the $\infty$-categories of derived stacks with values in $\infty$-categories)

$$\xymatrix{\underline{\mathcal{Q}\mathsf{Mod}} (n)\ar[rr] \ar[d] & & (\lqcoh) ^{\Delta[1]} \ar[d]^-{\textrm{ev}_0 \times \textrm{ev}_1} \\ \lqcoh \times \Delta[0] \ar[rr]_-{\mathrm{Sym}^{2} \times A[n]} & & \lqcoh \times \lqcoh  }$$

$$\xymatrix{\underline{\mathcal{Q}\mathsf{Perf}} (n)\ar[rr] \ar[d] & & (\lperf)^{\Delta[1]} \ar[d]^-{\textrm{ev}_0 \times \textrm{ev}_1} \\ \lperf \times \Delta[0] \ar[rr]_-{\mathrm{Sym}^{2} \times A[n]} & & \lperf \times \lperf  }$$ and analogous ones for stacks of nondegenerate objects.\\

By composing with the underlying space (or maximal $\infty$-subgroupoid) functor $\mathcal{I}$ (\cite[\S 1]{chern}), we obtain usual (i.e. with values in the $\infty$-category $\mathcal{S}$ of simplicial sets) derived stacks  
$$\mathcal{Q}\mathsf{Mod} (n)\, := \, \mathcal{I}(\underline{\mathcal{Q}\mathsf{Mod}} (n))\, \in \, \dst \, ,$$
$$\mathcal{Q}\mathsf{Mod}^{\mathrm{nd}} (n)\, := \, \mathcal{I}(\underline{\mathcal{Q}\mathsf{Mod}}^{\mathrm{nd}} (n))\, \in \, \dst \, ,$$
$$\mathcal{Q}\mathsf{Perf} (n)\, := \, \mathcal{I}(\underline{\mathcal{Q}\mathsf{Perf}} (n))\, \in \, \dst \, ,$$
$$\mathcal{Q}\mathsf{Perf}^{\mathrm{nd}} (n)\, := \, \mathcal{I}(\underline{\mathcal{Q}\mathsf{Perf}}^{\mathrm{nd}} (n))\, \in \, \dst \, .$$

If $X$ is a derived stack over $k$, let $\lqcoh(X)$ (resp. $\lperf (X)$) be the symmetric monoidal $\infty$-category of perfect (respectively, of quasi-coherent) complexes on $X$:
$$\lqcoh (X) := \lim_{\mathbb{R}\mathsf{Spec} (A) \to X} A-\idgmod$$
(resp. $$\lperf (X) := \lim_{\mathbb{R}\mathsf{Spec} (A) \to X} A-\iperf \,\,\, )$$
where both limits are taken in the $\infty$-category of stable symmetric monoidal presentable $\infty$-categories (see \cite{seattle}).
The $\infty$-functors $\bigwedge^2_A$, $\otimes_{A}$, $\mathrm{Sym}_{A}^2$ extend, by limits, to $\infty$-functors $\bigwedge^2_{\mathcal{O}_X} : \lqcoh (X) \longrightarrow \lqcoh (X)$, $\otimes_{{\mathcal{O}_X}}: \lqcoh (X) \times \lqcoh (X) \longrightarrow \lqcoh (X)$, $\mathrm{Sym}_{{\mathcal{O}_X}}^2: \lqcoh (X) \longrightarrow \lqcoh (X)$, and the full $\infty$-subcategory $\lperf(X)\subset \lqcoh(X)$, is stable under these $\infty$-functors.

\begin{df} \label{dbig}
Let $n\in \mathbb{Z}$, and $X$ be a derived stack over $k$. 
\begin{itemize}
\item The \emph{$\infty$-category of derived $n$-shifted quadratic complexes} over $X$ is given by the following pullback
$$\xymatrix{\underline{\mathsf{QMod}}(X;n)\ar[rr] \ar[d] & & \lqcoh (X)^{\Delta[1]} \ar[d]^-{\textrm{ev}_0 \times \textrm{ev}_1} \\ \lqcoh (X) \times \Delta[0] \ar[rr]_-{\mathrm{Sym}_{\mathcal{O}_X}^{2} \times \mathcal{O}_{X}[n]} & & \lqcoh (X) \times \lqcoh (X)  }$$ 
\item The \emph{$\infty$-category of derived $n$-shifted perfect quadratic complexes} over $X$ is given by the following pullback
$$\xymatrix{\underline{\mathsf{QPerf}}(X;n)\ar[rr] \ar[d] & & \lqcoh (X)^{\Delta[1]} \ar[d]^-{\textrm{ev}_0 \times \textrm{ev}_1} \\ \lperf (X) \times \Delta[0] \ar[rr]_-{\mathrm{Sym}_{\mathcal{O}_X}^{2} \times \mathcal{O}_{X}[n]} & & \lperf (X) \times \lperf (X)  }$$ 
\item The \emph{stack of derived $n$-shifted quadratic complexes} over $X$ is the object in $\dst$ $$\mathcal{Q}\mathsf{Mod}(X;n):=\mathrm{MAP}_{\dst}(X, \mathcal{Q}\mathsf{Mod}(n)).$$ 
\item The \emph{stack of derived $n$-shifted perfect quadratic complexes} over $X$ is the object in $\dst$ $$\mathcal{Q}\mathsf{Perf}(X;n):=\mathrm{MAP}_{\dst}(X, \mathcal{Q}\mathsf{Perf}(n)).$$
\item The \emph{classifying space of derived $n$-shifted quadratic complexes} over $X$ is the object in $\mathcal{S}$ $$\mathsf{QMod}(X;n):=\mathrm{Map}_{\dst}(X, \mathcal{Q}\mathsf{Mod} (n)).$$ 
\item The \emph{classifying space of derived $n$-shifted quadratic complexes} over $X$ is the object in $\mathcal{S}$ $$\mathsf{QPerf}(X;n):=\mathrm{Map}_{\dst}(X, \mathcal{Q}\mathsf{Perf} (n)).$$ 
\item The \emph{classifying space of derived $n$-shifted non degenerate quadratic complexes} over $X$ is the object in $\mathcal{S}$ $$\mathsf{QMod}^{\mathrm{nd}}(X;n):=\mathrm{Map}_{\dst}(X, \mathcal{Q}\mathsf{Mod}^{\mathrm{nd}} (n)).$$
\item The \emph{classifying space of derived $n$-shifted non degenerate quadratic complexes} over $X$ is the object in $\mathcal{S}$ $$\mathsf{QPerf}^{\mathrm{nd}}(X;n):=\mathrm{Map}_{\dst}(X, \mathcal{Q}\mathsf{Perf}^{\mathrm{nd}} (n)).$$
\end{itemize}
 \end{df}
 
 \begin{rmk}\label{harmony}\emph{Let $\mathrm{MAP}_{\dst^{\infty}}$ denote the internal Hom in the $\infty$-category $\dst^{\infty}$ of derived stacks (for the \'etale topology) with values in $\infty$-categories. We may view a usual derived stack $X$ (for the \'etale topology, and taking values in $\mathcal{S}$) as an object in $\dst^{\infty}$, by identifying $\mathcal{S}$ with the subcategory of $\infty$-groupoids inside the $\infty$-category of $\infty$-categories.  Since limits obviously commute with limits, and $\mathrm{MAP}_{\dst^{\infty}}(X, -)$ preserves fiber products, we have equivalences of $\infty$-categories $$\mathrm{MAP}_{\dst^{\infty}}(X,\underline{\mathcal{Q}\mathsf{Mod}} (n) ) (k) \simeq \underline{\mathsf{QMod}}(X;n)\, ,$$ $$\mathrm{MAP}_{\dst^{\infty}}(X,\underline{\mathcal{Q}\mathsf{Perf}} (n) ) (k) \simeq \underline{\mathsf{QPerf}}(X;n)\, .$$ Moreover, if $\mathcal{I}$ denote the maximal $\infty$-subgroupoid functor, we have equivalences in $\mathcal{S}$ $$\mathcal{I}(\underline{\mathsf{QMod}}(X;n)) \simeq \mathcal{Q}\mathsf{Mod}(X;n) (k) \simeq \mathsf{QMod}(X;n) \, ,$$ $$\mathcal{I}(\underline{\mathsf{QPerf}}(X;n)) \simeq \mathcal{Q}\mathsf{Perf}(X;n) (k) \simeq \mathsf{QPerf}(X;n).$$  }
 \end{rmk}


\bigskip

Consider the canonical projection $\infty$-functor $$\textrm{p}_1: \xymatrix{\underline{\mathsf{QMod}}(X;n) \ar[r] & \lqcoh (X)^{\Delta[1]}  \ar[r]^-{\textrm{ev}_0 } & \lqcoh (X) },$$ and, for $E \in \lqcoh (X)$ define $\infty$-category of \emph{$n$-shifted derived quadratic forms on} $E$ as the pullback $$\xymatrix{\underline{\mathbf{QF}}(E;n) \ar[r] \ar[d] & \underline{\mathsf{QMod}}(X;n) \ar[d]^-{\textrm{p}_1}\\ \Delta[0] \ar[r]_-{E} & \lqcoh (X). }$$ The underlying maximal $\infty$-subgroupoid $$\mathbf{QF}(E;n):=\mathcal{I}(\underline{\mathbf{QF}}(E;n))$$ is called the \emph{space of $n$-shifted derived quadratic forms on} $E$. Similar definitions for the $\infty$-category $\underline{\mathbf{QF}}^{\textrm{nd}}(E;n)$, and the space $\mathbf{QF}^{\textrm{nd}}(E;n)$, in the nondegenerate case.
Using the first point in Definition \ref{dbig} and Remark \ref{harmony}, we can make these more explicit as follows:


\begin{df}\label{global}Let $X$ be a derived stack over $k$, $E \in \lqcoh (X)$, and $n\in \mathbb{Z}$.
\begin{itemize}
\item The \emph{space of $n$-shifted derived quadratic forms on} $E$ is the mapping space $$\mathbf{QF}(E;n):= \mathrm{Map}_{ \lqcoh (X)}(\mathrm{Sym}_{\mathcal{O}_{X}}^2(E), \mathcal{O}_{X}[n])\,\,;$$
\item The \emph{set of $n$-shifted derived quadratic forms on} $E$ is $$\mathrm{QF}(E;n):= \pi_{0} \mathbf{QF}(E;n) $$ of connected components of $\mathbf{QF}(E;n)$, and an \emph{$n$-shifted quadratic form on} $E$ is by definition an element $q\in \mathrm{QF}(E;n)$.
\item The \emph{space} $\mathbf{QF}^{\textrm{nd}}(E;n)$ \emph{of $n$-shifted derived non-degenerate quadratic forms on} $E$ is defined by the following homotopy pullback diagram of simplicial sets
$$\xymatrix{\mathbf{QF}^{\textrm{nd}}(E;n) \ar[d] \ar[r] & \mathbf{QF}(E;n) \ar[d] \\ [\,\mathrm{Sym}_{\mathcal{O}_{X}}^2(E), \mathcal{O}_{X}[n]\,]^{\textrm{nd}} \ar[r] & [\, \mathrm{Sym}_{\mathcal{O}_{X}}^2(E), \mathcal{O}_{X}[n]\,] }$$ where $[-,-]$ denotes the hom-sets in the homotopy category of $\lqcoh (X)$, and $ [\,\mathrm{Sym}_{\mathcal{O}_{X}}^2(E), \mathcal{O}_{X}[n]\,]^{\textrm{nd}}$ is the subset of $ [\,\mathrm{Sym}_{\mathcal{O}_{X}}^2(E), \mathcal{O}_{X}[n]\,]$ consisting of maps $v: \mathrm{Sym}_{\mathcal{O}_{X}}^2(E) \to \mathcal{O}_{X}[n]$ such that the adjoint map $v^{\flat}: E \to E^{\vee}[n]$ is an isomorphism in $\ho (\lqcoh (X))$.
\item The \emph{set} $\mathrm{QF}^{\textrm{nd}}(E;n)$ \emph{of $n$-shifted derived non-degenerate quadratic forms on} $E$ is the set $$\mathrm{QF}^{\textrm{nd}}(E;n):= \pi_{0} \mathbf{QF}(E;n)^{\textrm{nd}} $$ of connected components of $\mathbf{QF}(E;n)^{\textrm{nd}}$, and an \emph{$n$-shifted non-degenerate quadratic form on} $E$ is by definition an element $q\in \mathrm{QF}(E;n)^{\textrm{nd}}$.
\item an $n$\emph{-shifted derived} (resp. \emph{non-degenerate}) \emph{quadratic complex on} $X$ is a pair $(E,q)$ where $E\in \lqcoh (X)$ and $q\in \mathrm{QF}(E;n)$ (resp. $q\in \mathrm{QF}^{\textrm{nd}}(E;n)$).

\end{itemize} 
\end{df}

 \begin{rmk}\emph{Though Definition \ref{global} makes sense for any $E \in \lqcoh (X)$, we will be mostly interested in the case where $E \in \lperf (X)$. }
\end{rmk}

Note that, by definition, a derived $n$-shifted quadratic form $q$ on $E$ is a map $q: \mathrm{Sym}_{\mathcal{O}_{X}}^2(E) \to \mathcal{O}_{X}[n]$ in the homotopy category $\textrm{h}(\lqcoh (X))\simeq\mathrm{D}_{\mathrm{QCoh}}(X)$ (equivalence of triangulated $k$-linear categories, since the $\infty$-category $\lqcoh (X)$ is $k$-linear and stable). \\

\begin{ex}\emph{A \emph{symmetric obstruction theory}, according to \cite[Def. 1.10]{befa}, is an example of a derived $1$-shifted quadratic complex. }
\end{ex}

\begin{df}\label{resglobal} If $E_1 \in \lqcoh (X)$, $(E_2,q_2)$ is a derived $n$-shifted quadratic complex over $X$, and $f:E_1 \to E_2$ is a map in $\lqcoh (X)$, then the composite $$\xymatrix{\mathrm{Sym}_{\mathcal{O}_{X}}^{2}E_1 \ar[rr]^-{\mathrm{Sym}^{2} f} & & \mathrm{Sym}_{\mathcal{O}_{X}}^{2}E_2 \ar[r]^-{q_2} & \mathcal{O}_{X}[n]}$$ defines an $n$-shifted quadratic form on $E_1$, that we denote by $f^*q_2$. $f^*q_2$ is called the \emph{pull-back} or \emph{restriction} of $q_2$ along $f$.
\end{df}

The next result establishes a first link between derived symplectic structures and (non-degenerate) derived quadratic forms on derived stacks. We use the notations of \cite{ptvv}.\\
 
\begin{prop} \label{sympltoquad} Let be a derived lfp stack over $k$, $\mathrm{char} \, k =0$. There are canonical maps in $\ho (\ssets)$ $$\mathbf{Sympl}(X;n) \longrightarrow \mathbf{QF}^{\textrm{nd}}(\mathbb{T}_{X}[1];n-2), $$ $$\mathbf{Sympl}(X;n) \longrightarrow \mathbf{QF}^{\textrm{nd}}(\mathbb{T}_{X}[-1];n+2). $$

\end{prop}

\noindent \textbf{Proof.} It's enough to recall that $\wedge_{\mathcal{O}_{X}}^{2} \mathbb{T}_{X} \simeq \mathrm{Sym}^{2}_{\mathcal{O}_{X}}(\mathbb{T}_{X}[\pm1])[\mp2]$; this yields maps (actually equivalences in $\mathcal{S}$) $$\mathcal{A}^{2,\textrm{nd}}(X; n) \longrightarrow \mathbf{QF}^{\textrm{nd}}(\mathbb{T}_{X}[\pm1]; n\mp2),$$ and we just precompose each of these with the canonical underlying-$2$-form map $\mathbf{Sympl}(X;n) \longrightarrow \mathcal{A}^{2,\textrm{nd}}(X; n)$.

\hfill $\Box$

\begin{rmk}\label{moreprecisely}\emph{ Note that, more generally, for any $m \in \mathbb{Z}$, the d\'ecalage isomorphisms $\wedge_{\mathcal{O}_{X}}^{2} \mathbb{T}_{X} \simeq \mathrm{Sym}^{2}_{\mathcal{O}_{X}}(\mathbb{T}_{X}[2m+1])[-4m-2]$ yield equivalences in $\mathcal{S}$ $$\mathcal{A}^{2}(X; n) \simeq \mathbf{QF}(\mathbb{T}_{X}[2m+1]; n-4m-2)$$ between $n$-shifted $2$-forms on $X$ (\cite[Def. 1.12 and Prop. 1.14]{ptvv}) and $(n-4m-2)$-shifted quadratic forms on $\mathbb{T}_{X}[2m+1]$. }

\end{rmk}

As done in \cite[Def. 1.10]{ptvv} for shifted symplectic forms, we may give the following
 
  \begin{df} For a derived $m$-shifted quadratic (not necessarily non-degenerate) form $q$ on $\mathbb{T}_{X}[\pm 1]$, we define the \emph{space of keys} of $q$, as the homotopy fiber at $q$ of the composite map $$\mathcal{A}^{2, \textrm{cl}}(X; m\pm 2) \longrightarrow \mathcal{A}^{2}(X; m\pm 2) \longrightarrow \mathbf{QF}(\mathbb{T}_{X}[\pm1]; m).$$
\end{df}

Exactly in the same way as done in  Section \ref{dercliff} for the case of complexes over simplicial commutative algebra (i.e the case $X=\mathbb{R}\mathsf{Spec} (A)$), for a map of derived stacks $\varphi: Y\to X$, and $E\in \lqcoh (X)$   we may define the \emph{base-change} maps $\ho (\ssets)$ $$\varphi^*: \mathbf{QF}_{X}(E;n)\longrightarrow \mathbf{QF}_{Y}(\varphi^{*}E;n),$$ and 
$$\varphi^*: \mathbf{QF}_{X}^{\textrm{nd}}(E;n)\longrightarrow \mathbf{QF}_{Y}^{\textrm{nd}}(\varphi^{*}E;n).$$\\

Similarly, we define the \emph{orthogonal sum} $\perp$ of derived $n$-shifted quadratic complexes over $X$ and the notion of \emph{derived isometric structure} on a map $f: (E,q)\to (E',q')$ between derived quadratic complexes in $\lqcoh (X)$, by globalizing Definition \ref{sum}  and \ref{d3}.

\subsection{Derived Grothendieck-Witt groups} \label{GW}

In this section we define a derived version of the Grothendieck-Witt group of a derived stack. Its functoriality and further applications will be given elsewhere.\\

Recall from Def. \ref{dbig} the classifying space $\mathsf{QPerf}(X;n):=\mathrm{Map}_{\dst}(X, \mathcal{Q}\mathsf{Perf}^{\mathrm{nd}} (n))$ (resp. $\mathsf{QPerf}^{\mathrm{nd}}(X;n):=\mathrm{Map}_{\dst}(X, \mathcal{Q}\mathsf{Perf}^{\mathrm{nd}} (n))$) of derived $n$-shifted (resp. non degenerate) quadratic complexes over a lfp derived stack $X$.
The orthogonal sum $\perp$ of derived quadratic complexes induces both on $\pi_{0}(\mathsf{QPerf}(X;n))$ and on $\pi_{0}(\mathsf{QPerf}^{\mathrm{nd}}(X;n))$ a commutative monoid structure, still denoted by $\perp$. 

\begin{df} 
Let $n\in \mathbb{Z}$, and $X$ be a derived Artin stack lfp over $k$. The \emph{extended derived Grothendieck-Witt group} of $X$ is the Grothendieck group of the commutative monoid $\pi_{0}(\mathsf{QPerf}(X;n))$
$$\widehat{\mathsf{GW}^{\mathrm{ext}}}(X;n):= K_0 (\pi_{0}(\mathsf{QC}(X;n)), \perp).$$ The \emph{derived Grothendieck-Witt group} of $X$ is the Grothendieck group of the commutative monoid $\pi_{0}(\mathsf{QPerf}^{\mathrm{nd}}(X;n))$
$$\widehat{\mathsf{GW}}(X;n):= K_0 (\pi_{0}(\mathsf{QPerf}(X;n)^{\mathrm{nd}}), \perp).$$
 \end{df}
 
 \begin{rmk}\emph{Unlike the classical, unshifted case, the derived (extended) Grothendieck-Witt will only be a ring if we consider all (or just all \emph{even}) shifts at the same time, the tensor product of an $n$-shifted quadratic complex with an $m$-shifted quadratic complex being naturally an $(n+m)$-shifted quadratic complex.}
 \end{rmk}

When $X= \mathrm{Spec} \, k$ and $n=0$, one easily verifies that $(\pi_{0}(\mathsf{QPerf}^{\mathrm{nd}}(X;n)), \perp)$ is isomorphic (as a commutative monoid) to the classical monoid of isomorphism classes of finitely generated projective non-degenerate quadratic $k$-modules under orthogonal sum (see, e.g. \cite[\S 1.8]{mire}). Therefore, we get 

\begin{prop}\label{gwcomp} We have a canonical isomorphism of abelian groups  $$\widehat{\mathsf{GW}}(k;0) \simeq \widehat{GW}(k),$$ where $ \widehat{GW}(k)$ denotes the classical Grothendieck-Witt group of $k$.
\end{prop}

\subsection{Derived quadratic stacks }

\begin{df} Let $X$ be a derived Artin stack locally finitely presented ($\equiv$ lfp) over $k$, and $n\in \mathbb{Z}$.
\begin{itemize}
\item The \emph{space of $n$-shifted derived quadratic forms over} $X$ is the space $$\mathbf{QF}(X;n):=\mathbf{QF}(\mathbb{T}_{X};n).$$
\item The \emph{set of $n$-shifted derived quadratic forms over} $X$ is $$\mathrm{QF}(X;n):= \pi_{0} \mathbf{QF}(X;n) $$ of connected components of $\mathbf{QF}(X;n)$, and an \emph{$n$-shifted quadratic form over} $X$ is by definition an element $q\in \mathrm{QF}(X;n)$.
\item The space $\mathbf{QF}^{\textrm{nd}}(X;n):=\mathbf{QF}^{\textrm{nd}}(\mathbb{T}_{X};n) $ is the space \emph{of $n$-shifted derived non-degenerate quadratic forms over} $X$ 
\item The \emph{set} $\mathrm{QF}^{\textrm{nd}}(X;n):= \pi_{0} \mathbf{QF}(E;n)^{\textrm{nd}}$ is the set \emph{of $n$-shifted derived non-degenerate quadratic forms over} $X$,  and an \emph{$n$-shifted non-degenerate quadratic form over} $X$ is by definition an element $q\in \mathrm{QF}(X;n)^{\textrm{nd}}$.
\item an $n$\emph{-shifted derived} (resp. \emph{non-degenerate}) \emph{quadratic stack} is a pair $(X,q)$, where $X$ is a derived stack locally of finite presentation over $k$, and $q\in \mathrm{QF}(X;n)$ (resp. $q\in \mathrm{QF}^{\textrm{nd}}(X;n)$).

\end{itemize} 
\end{df}

\begin{rmk}\emph{\textsf{Derived moduli stack of shifted symplectic structures.}  Given $n \in \mathbb{Z}$ a derived lfp stack  $X$ over $k$, $\textrm{char}\, k =0$, one can consider the derived moduli stack $\mathsf{Sympl}(X; n)$ of $n$-shifted symplectic structures on $X$, as the functor sending $A \in \salg$ to the space of $n$-shifted, relative symplectic structures on $X \times \mathbb{R}\mathrm{Spec} \, A \to \mathbb{R}\mathrm{Spec} \, A$. The derived stack $\mathsf{Sympl}(X; n)$ can be shown to have a tangent complex $\mathbb{T}_{\mathsf{Sympl}(X; n)}$ such that $\mathbb{T}_{\mathsf{Sympl}(X; n), \, \omega} \simeq \mathcal{A}^{2,\, \textrm{cl}}(X)[n]$ (the complex of $k$-dg modules of  $n$-shifted closed $2$-forms on $X$), at any $k$-point (i.e. $n$-shifted symplectic structure on $X$) $\omega$. I conjecture that, under suitable hypotheses on $X$, one can generalize the classical $C^{\infty}$ arguments of \cite{fricke} in order to show that $\mathsf{Sympl}(X; n)$ is in fact a shifted quadratic derived stack. }
\end{rmk}

\begin{df}\label{shriek} If $X_1$ is a derived Artin stack lfp over $k$, $(X_2,q_2)$ an $n$-shifted derived quadratic stack over $k$, and $f:X_1 \to X_2$ is a map in $\ho (\dst)$, then the composite 
$$\xymatrix{\mathrm{Sym}_{\mathcal{O}_{X_1}}^ {2} \mathbb{T}_{X_1} \ar[r]  & \mathrm{Sym}_{\mathcal{O}_{X_1}}^{2}f^* \mathbb{T}_{X_2} \ar[r]^-{\sim} & f^*\mathrm{Sym}_{\mathcal{O}_{X_2}}^{2}\mathbb{T}_{X_2}  \ar[r]^-{f^* q_2} & f^* \mathcal{O}_{X_2}[n] \simeq \mathcal{O}_{X_1}[n]}$$ defines an $n$-shifted quadratic form on $X_1$, that we denote by $f^{q}q_2$.  More generally, we still denote by $f^{q}$ the induced map $$f^{q}: \mathbf{QF}(X_2;n) \longrightarrow \mathbf
{QF}(X_1;n)$$ in $\ho (\ssets)$.
\end{df}

\begin{rmk} \emph{Note that $f^{q} q_2$ has the following equivalent description. Let us denote by $\mathrm{T}f: \mathbb{T}_{X} \to f^{*}\mathbb{T}_{Y}$ the induced tangent map. Then $$f^{q} q= (\mathrm{T}f)^{*} (f^* q)$$  were $f^* q$ denotes the base-change of $q$ (a quadratic form on $f^* \mathbb{T}_{Y}$ over $X$), and $(\mathrm{T}f)^{*}(\varphi^* q)$ the restriction of $f^* q$ along $f$ (as in Definition \ref{resglobal}).}
\end{rmk}

The derived versions of $f$ being an isometry or a null-map, are not properties of but rather data on $f$. More precisely, we give the following

\begin{df} Let $(X_1, q_1)$ and  $(X_2,q_2)$ be $n$-shifted derived quadratic stack sover $k$ \begin{itemize}
\item The \emph{space of derived isometric structures} on a map $f:X_1\to X_2$ is by definition the space $$\mathsf{Isom}(f;(X_1,q_1),(X_2,q_2)):= \mathrm{Path}_{q_1, f^*q_2}(\mathbf{QF}(X_1;n)).$$ 
\item The space $\mathsf{Null}(f; X_1, X_2, q_2):=\mathsf{Isom}(f;(X_1,0),(X_2,q_2))$ is called the \emph{space of derived null-structures} on $f$.
\item A \emph{derived isometric structure} on a map $f:X_1\to X_2$ is an element in $\pi_{0}(\mathsf{Isom}(f;(X_1,q_1),(X_2,q_2)))$. 
\item A \emph{derived null structure} on a map $f:X_1\to X_2$ is an element in $\pi_{0}(\mathsf{Isom}(f;(X_1,0),(X_2,q_2)))$. 
\end{itemize}
\end{df}

Now we want to define a condition of \emph{non-degeneracy} on a null-structure on a given map $f:X_1 \to (X_2,q_2)$ that will prove useful later. This is a quadratic analog of the notion of Lagrangian structure for derived symplectic forms (\cite[\S 2.2]{ptvv}).
Let $\gamma \in \mathsf{Null}(f; X_1, X_2, q_2)$ be a fixed null-structure on $f$.
By definition, $\gamma$ is a path between $0$ and the composite morphism
$$\xymatrix{\mathrm{Sym}_{\mathcal{O}_{X_1}}^2 \mathbb{T}_{X_1} \ar[r] & 
f^{*} (\mathrm{Sym}_{\mathcal{O}_{X_2}}^2\mathbb{T}_{X_2})  \ar[r] & \mathcal{O}_{X_1}[n].}$$
If $\mathbb{T}_{f}$ is the relative tangent complex of $f$, so that we
 have the transitivity exact triangle of perfect complexes on $X_1$
$$\mathbb{T}_{f} \longrightarrow \mathbb{T}_{X_1} \longrightarrow f^{*}(\mathbb{T}_{X_2}).$$
The null-structure $\gamma$ induces also a path $\gamma'$ between $0$ and the composite morphism
$$\xymatrix{\varphi:\mathbb{T}_{f} \otimes \mathbb{T}_{X_1} \ar[r] & \mathbb{T}_{X_1} \otimes \mathbb{T}_{X_1} \ar[r] & \mathrm{Sym}_{\mathcal{O}_{X_1}}^2 \mathbb{T}_{X_1}  \ar[r]^-{f^* q_2} & \mathcal{O}_{X_1}[n].}$$
But the morphism
$\mathbb{T}_{f} \longrightarrow f^{*}(\mathbb{T}_{X_2})$ comes itself with 
a canonical (independent of $\gamma$) path from itself to $0$, so we get another induced path $\delta$ from $\varphi$ to $0$.
By composing $\gamma'$ and $\delta$, we then obtain a loop pointed at $0$ in the space 
$$\mathrm{Map}_{\lqcoh(X_1)}(\mathbb{T}_{f} \otimes \mathbb{T}_{X_1},\mathcal{O}_{X_1}[n]).$$ This loop defines
an element in 
$$\pi_{1}(\mathrm{Map}_{\lqcoh(X_1)}(\mathbb{T}_{f} \otimes \mathbb{T}_{X_1},\mathcal{O}_{X_1}[n]) ; 0) \simeq
\mathrm{Hom}_{\mathsf{D}_{\mathrm{QCoh}}}(\mathbb{T}_{f} \otimes \mathbb{T}_{X_1},\mathcal{O}_{X_1}[n-1]].$$
By adjunction, we get a morphism of perfect complexes on $X_1$
$$\Theta_{\gamma} : \mathbb{T}_{f} \longrightarrow \mathbb{L}_{X_1}[n-1],$$
depending on the null-structure structure $\gamma$. 

\begin{df}\label{ndnull}
Let $f : X_1 \longrightarrow X_2$ be a morphism of derived Artin stacks and $q_2$ a derived $n$-shifted quadratic form on $X_2$. An null-structure $\gamma$ on $f: X_1 \to (X_2,q_2)$ is \emph{non-degenerate} if
the induced morphism
$$\Theta_{\gamma} : \mathbb{T}_{f} \longrightarrow \mathbb{L}_{X_1}[n-1]$$
is an isomorphism in $\mathsf{D}(X_1)$.
\end{df}

\section{Existence theorems}

In this Section we prove two existence theorems for derived quadratic forms on derived stacks, directly inspired by the analogous results for derived symplectic structures (\cite[Thms. 2.5 and 2.9]{ptvv}). For notation and definitions of $\mathcal{O}$-compact derived stack and $\mathcal{O}$-orientation, we refer the reader to \cite[\S 2]{ptvv}.  \\

\begin{thm}\label{map} Let $X$ be a derived stack locally of finite presentation over $k$, $(m,n)\in \mathbb{Z}^2$, and $q \in QF(\mathbb{T}_{X}[m];n)$. Let $Y$ be an $\mathcal{O}$-compact derived stack equipped with an $\mathcal{O}$-orientation $\eta$ of degree $d$ , and suppose that the derived mapping stack $\textsf{MAP}_{\dst}(Y,X)$ is a derived Artin stack locally of finite presentation over $k$. Then,  the $m$-shifted tangent complex $\mathbb{T}_{\textsf{MAP}_{\dst}(Y,X)}[m]$ admits a canonical $(n-d)$-shifted derived quadratic form $q'\equiv q'_{\eta}$. If moreover $q$ is non-degenerate, then so is $q'_{\eta}$.
\end{thm}


\textbf{Proof.} We borrow the notations from \cite[\S 2.1]{ptvv}, and define $q'$ as follows.
For $x : \spec\, A \longrightarrow \mathsf{Map}(Y,X)$ an $A$-point corresponding to a morphism of derived 
stacks
$$f : Y_{A}:=Y\times \spec\, A \longrightarrow X\, ,$$
the tangent complex of $\textsf{MAP}_{\dst}(Y,X)$ at the point $x$ is given by 
$$\mathbb{T}_{x}\textsf{MAP}_{\dst}(Y,X) \simeq \mathbb{R}\mathrm{Hom}(\mathcal{O}_{Y_{A}},f^{*}(\mathbb{T}_{X})).$$
The quadratic form  
$$q: \mathrm{Sym}^{2}_{\mathcal{O}_{X}}(\mathbb{T}_{X}[m])  \longrightarrow \mathcal{O}_{X}[n],$$ 
induces by pullback a map of $A$-dg-modules
$$\rho_{A}: \mathrm{Sym}^{2}_{A}(\mathbb{R}\mathrm{Hom}(\mathcal{O}_{Y_{A}},f^{*}(\mathbb{T}_{X}))[m]) \longrightarrow \mathbb{R}\mathrm{Hom}(\mathcal{O}_{Y_{A}},\mathcal{O}_{Y_{A}}[n]) .$$
Now we use the $d$-orientation $\eta$ on $Y$ (\cite[Def. 2.4]{ptvv}): its $[n]$-shift induces, by definition a map of $A$-dg modules
$$\eta_{A}[n] : \mathbb{R}\mathrm{Hom}(\mathcal{O}_{Y_{A}},\mathcal{O}_{Y_{A}}[n]) \longrightarrow A[n-d]$$

By composing $\eta_{A}[n]$ with $\rho_{A}$, we get 
$$q'_{A}: \mathrm{Sym}^{2}_{A} (\mathbb{T}_{x}\textsf{MAP}_{\dst}(Y,X)[m]) \simeq \mathrm{Sym}^{2}_{A}(\mathbb{R}\mathrm{Hom}(\mathcal{O}_{Y_{A}},f^{*}(\mathbb{T}_{X}))[m]) \longrightarrow A[n-d].$$
This defines the quadratic form $q'\equiv q'_{\eta}$ on $\mathbb{T}_{\textsf{MAP}_{\dst}(Y,X)}[m]$. 
When $q$ is non-degenerate, then $q'_{\eta}$ is also non-degenerate, since $\eta$ is an orientation.
\hfill $\Box$ \\

\begin{prop}\label{perf&BG}
Let $\mathrm{B}GL_{n}$ be the classifying stack of $GL_{n,k}$, and $\mathbb{R}\mathbf{Perf}$ the derived stack of perfect complexes (sending a simplicial $k$-algebra $A$ to the nerve of the category of cofibrant perfect $A$-dg-modules with equivalences between them; see \cite[\S 2.3]{ptvv}). Then, for any $m\in \mathbb{Z}$, $\mathbb{T}_{\mathrm{B}GL_{n}}[2m+1]$, and $\mathbb{T}_{\mathbb{R}\mathbf{Perf}}[2m+1]$ are equipped with canonical  derived $(2(2m+2))$-shifted quadratic forms.
\end{prop}

\textbf{Proof.} This follows from the end of \cite[\S 1.2]{ptvv} (for $\mathrm{B}GL_n$), and from \cite[Thm. 2.12]{ptvv} (for $\mathbb{R}\mathbf{Perf}$), via our Proposition \ref{sympltoquad} and Remark \ref{moreprecisely}. \\Alternatively, one can easily and more directly verify that
\begin{itemize}
\item $(A,B) \mapsto \mathsf{tr}(AB)$ defines a non-degenerate quadratic form $\mathrm{Sym}^2(\mathsf{Lie}(GL_{n,k})) \to k$ on $\mathsf{Lie}(GL_{n,k})$, and 
\item if $\mathcal{E}$ is the universal perfect complex on $\mathbb{R}\mathbf{Perf}$, and $\mathcal{A}:= \mathbb{R}\underline{\mathrm{End}}(\mathcal{E}) \simeq \mathcal{E} \otimes \mathcal{E}^{\vee}$, then the composite
$$\xymatrix{\mathcal{A} \otimes \mathcal{A} \ar[r]^-{\mu} & \mathcal{A}\simeq \mathcal{E} \otimes \mathcal{E}^{\vee} \ar[r]^-{\mathrm{ev}} & \mathcal{O}_{\mathbb{R}\mathbf{Perf}} }$$ (where $\mu$ is the multiplication and $ev$ is the canonical pairing), induces a non-degenerate derived quadratic form $\mathrm{Sym}_{\mathcal{O}_{\mathbb{R}\mathbf{Perf}}}^2 (\mathcal{A}) \to \mathcal{O}_{\mathbb{R}\mathbf{Perf}}$ on $\mathcal{A}$.
\end{itemize}
Now, the Proposition follows from the decalage isomorphisms $\mathrm{Sym}^2(\mathbb{T}[2m+1]) \simeq \wedge^{2}(\mathbb{T})[2(2m+1)]$, and from the identifications (see \cite{ptvv}) $$\mathbb{T}_{\mathrm{B}GL_n} \simeq \mathsf{Lie}(GL_{n,k})[1]$$ $$\mathbb{T}_{\mathbb{R}\mathbf{Perf}} \simeq \mathcal{A}[1].$$
 \hfill $\Box$ \\
 
 \begin{rmk} \emph{Proposition \ref{perf&BG} remains true by replacing $GL_{n}$ with any affine smooth reductive group scheme $G$ over $k$: it is enough to use any non-degenerate $G$-invariant bilinear form on the Lie algebra $\mathsf{Lie}(G)$. See the end of \cite[\S1.2]{ptvv} for details. }
 \end{rmk}

\begin{cor}\label{mapcor} \begin{enumerate} \item Let $(X,q)$ be an $n$-shifted quadratic derived stack, $Y$ be an $\mathcal{O}$-compact derived stack equipped with an $\mathcal{O}$-orientation $\eta$ of degree $d$ , and suppose that the derived mapping stack $\textsf{MAP}_{\dst}(Y,X)$ is a derived Artin stack locally of finite presentation over $k$. Then,  $\textsf{MAP}_{\dst}(Y,X)$ admits a canonical $(n-d)$-shifted derived quadratic form $q'\equiv q'_{\eta}$. If moreover $q$ is non-degenerate, then so is $q'_{\eta}$.
\item In the following cases, there exists a canonical derived $n$-shifted quadratic form on $\mathbb{T}_{Z}[2m+1]$ (for the De Rham and Dolbeault notations below, see e.g. \cite{simp} or \cite{ptvv}):
\begin{itemize}
\item $Z:=\mathbb{R}\mathbf{Perf}_{DR}(Y):=\textsf{MAP}_{\dst}(Y_{DR},\mathbb{R}\mathbf{Perf})$ is the derived stack 
of perfect complexes with flat connections on $Y$, $m \in \mathbb{Z}$, $n:= 2(2m+2-d)$, where $Y$ is a smooth and proper
Deligne-Mumford stack with connected geometric fibers of relative dimension $d$, together with a choice of a fundamental de Rham class $\eta_{DR, Y}\equiv [Y]_{DR} \in H^{2d}_{DR}(Y,\mathcal{O})$;
\item $Z:=\mathbb{R}\mathbf{Perf}_{\textrm{Dol}}(Y):=\textsf{MAP}_{\dst}(Y_{\textrm{Dol}},\mathbb{R}\mathbf{Perf})$ is the derived stack 
of perfect complexes with Higgs fields on $Y$, $m \in \mathbb{Z}$, $n:= 2(2m+2-d)$, where $Y$ is a smooth and proper
Deligne-Mumford stack with connected geometric fibers of relative dimension $d$, together with a choice of a fundamental Dolbeault class $\eta_{\textrm{Dol}, Y}\equiv [Y]_{\textrm{Dol}} \in H^{2d}_{\textrm{Dol}}(Y,\mathcal{O})$;
\item $Z:=\mathbb{R}\mathbf{Perf}(Y):=\textsf{MAP}_{\dst}(Y,\mathbb{R}\mathbf{Perf})$ is the derived stack 
of perfect complexes on $Y$, $m \in \mathbb{Z}$, $n:= 2(2m+2)-d$, where $Y$ is a smooth and proper Calabi-Yau
Deligne-Mumford stack with connected geometric fibers of relative dimension $d$, together with a choice of a trivialization of the canonical sheaf $\eta_{Y}: \omega_{Y/k} \simeq \mathcal{O}_{Y}$;
\item  $Z:=\mathbb{R}\mathbf{Perf}(M):=\textsf{MAP}_{\dst}(\mathrm{Sing(M)},\mathbb{R}\mathbf{Perf})$ is the derived stack 
of perfect complexes on $M$, $m \in \mathbb{Z}$, $n:= 2(2m+2)-d$, where $M$ is a compact oriented topological manifold of dimension $d$, together with a choice of a topological fundamental class $[M]\in H_{d}(M,k)$;

\item $Z:=\mathbb{R}\mathbf{Vect}_{r, DR}(Y):=\textsf{MAP}_{\dst}(Y_{DR},\mathrm{B}GL_r)$ is the derived stack 
of rank $r$ vector bundles  with flat connections on $Y$, $m \in \mathbb{Z}$, $n:= 2(2m+2-d)$, where $Y$ is a smooth and proper
Deligne-Mumford stack with connected geometric fibers of relative dimension $d$, together with a choice of a fundamental de Rham class $\eta_{DR, Y}\equiv [Y]_{DR} \in H^{2d}_{DR}(Y,\mathcal{O})$;
\item $Z:=\mathbb{R}\mathbf{Vect}_{r, \textrm{Dol}}(Y):=\textsf{MAP}_{\dst}(Y_{\textrm{Dol}}, \mathrm{B}GL_r)$ is the derived stack 
of rank $r$ vector bundles with Higgs fields on $Y$, $m \in \mathbb{Z}$, $n:= 2(2m+2-d)$, where $Y$ is a smooth and proper
Deligne-Mumford stack with connected geometric fibers of relative dimension $d$, together with a choice of a fundamental Dolbeault class $\eta_{\textrm{Dol}, Y}\equiv [Y]_{\textrm{Dol}} \in H^{2d}_{\textrm{Dol}}(Y,\mathcal{O})$;
\item $Z:=\mathbb{R}\mathbf{Vect}_{r}(Y):=\textsf{MAP}_{\dst}(Y,\mathrm{B}GL_r)$ is the derived stack 
of rank $r$ vector bundles on $Y$, $m \in \mathbb{Z}$, $n:= 2(2m+2)-d$, where $Y$ is a smooth and proper Calabi-Yau
Deligne-Mumford stack with connected geometric fibers of relative dimension $d$, together with a choice of a trivialization of the canonical sheaf $\eta_{Y}: \omega_{Y/k} \simeq \mathcal{O}_{Y}$;
\item  $Z:=\mathbb{R}\mathbf{Loc}_{r}(M):=\textsf{MAP}_{\dst}(\mathrm{Sing(M)},\mathrm{B}GL_r)$ is the derived stack 
of rank $r$ local systems of $k$-vector spaces on $M$, $m \in \mathbb{Z}$, $n:= 2(2m+2)-d$, where $M$ is a compact oriented topological manifold of dimension $d$, together with a choice of a topological fundamental class $[M]\in H_{d}(M,k)$;
\end{itemize}
\end{enumerate}
\end{cor}

\textbf{Proof.} Part (1) is just Theorem \ref{map} in the case $m=0$. Part (2) follows from Proposition \ref{perf&BG} and Theorem \ref{map}, given that under the respective given hypotheses $Y_{DR}$ and $Y_{\textrm{Dol}}$ are $\mathcal{O}$-compact derived stack equipped with $\mathcal{O}$-orientations $[Y]_{DR}$ and $[Y]_{\textrm{Dol}}$ of degree $2d$, while $Y$ is $\mathcal{O}$-compact derived stack equipped, via Serre duality, with an $\mathcal{O}$-orientation $\eta$ of degree $d$, in the Calabi-Yau case (see \cite{ptvv}).
 \hfill $\Box$ \\

 \begin{rmk} \label{chiariamo}\emph{Part (2) of Corollary \ref{mapcor} also follows directly from the finer corresponding results of \cite[Cor. 2.13]{ptvv}, via the use of our Proposition \ref{sympltoquad} and Remark \ref{moreprecisely}. However the proof in our quadratic case is considerably simpler, since we are not interested in proving the existence of any closedness data on the associated shifted $2$-form.}
 \end{rmk}

\begin{cor}\label{loop} If $(X,q)$ is an $n$-shifted quadratic derived stack, its derived loop stack $\mathsf{L} X := \textsf{MAP}_{\dst}(S^1,X)$ has an induced derived $(n-1)$-shifted quadratic form.
\end{cor}

\textbf{Proof.} This follows from Corollary \ref{mapcor} (1), since for any compact, oriented $d$- dimensional manifold $M$, the constant derived stack with value its singular simplicial set $ \mathrm{Sing}(M)$ is canonically a $\mathcal{O}$-compact derived stack equipped with an $\mathcal{O}$-orientation $\eta$ of degree $d$ induced from the fundamental class. \\ Another proof can be given using Theorem \ref{lagr} below, by noticing that $$\mathsf{L} X\simeq X\times_{X\times X} X$$ and endowing $X \times X$ with the quadratic form $(q, -q)$. We leave the details of this alternative proof to the interested reader. 
 \hfill $\Box$ \\


\begin{thm}\label{lagr} Let $(X,q)$ be an $n$-shifted quadratic derived stack, and $(f_i: Y_i \to X, \gamma_i)$, $i=1,2$ two null-structures relative to $(X,q)$. Then, the homotopy fiber product $Y_1 \times^h_X Y_2$ admits a canonical $(n-1)$-shifted derived quadratic form. If moreover, the null structures $(f_i: Y_i \to X, \gamma_i)$, $i=1,2$ are non-degenerate, and $(X,q)$ is non-degenerate, so is $Y_1 \times^h_X Y_2$ with this induced $(n-1)$-shifted derived quadratic form.
\end{thm}

\noindent \textbf{Proof.} The proof consists in decoupling, in the proof of \cite[Thm. 2.9]{ptvv}, the existence of a $2$-form and its non-degneracy from the additional closedness data.
 To ease notations we will write, for the duration of the proof, $g^{*}$ instead of $g^{q}$ (Def. \ref{shriek}), for an arbitrary map $g$ in $\dst$. \\ Let $Z:=Y_1\times_{X}^{h}Y_2$.  By definition of homotopy fiber product, the two morphisms
$$p_1 : \xymatrix{Z \ar[r] & Y_1 \ar[r]^{f_1} & X} \qquad p_2 : \xymatrix{Z \ar[r] & Y_2 \ar[r]^-{f_2} & X}$$
come equipped with a natural path $u$between them. Now, $u$ gives rise to 
a path between the induced morphisms on the spaces of derived $n$-shifted quadratic forms
$$u^{*} : p_1^{*} \leadsto p_2^{*} : \mathbf{QF}(X;n) \longrightarrow 
\mathbf{QF}(Z;n).$$
Moreover, $\gamma_i$ defines a path in the space $\mathbf{QF}(Z;n)$ between $0$ and $p^{*}_i q$, for $i=1,2$.
By composing $\gamma_1$, $u^{*}(q)$ and $\gamma_{2}^{-1}$, we get 
a loop at $0$ in the space $\mathbf{QF}(Z;n)$, thus a well
defined element
$$Q=Q(q,\gamma_1,\gamma_2) \in \pi_{1}(\mathbf{QF}(Z;n);0)\simeq \pi_{0}(\mathbf{QF}(Z;n-1))= \mathrm{QF}(Z;n-1).$$
Let us now suppose that $(X,q)$ is non-degenerate and so are the null-structures.  We will prove that $Q= Q(q,\gamma_1,\gamma_2)$ is also non-degenerate.
Let $\Theta_{\gamma_i}:\mathbb{T}_{f_i} \longrightarrow \mathbb{L}_{Y_i}[n-1]$, $i=1,2$ be the induced maps (see Def. \ref{ndnull}), and $pr_{i} : Z \longrightarrow X_i$, $i=1,2$ the two projections.
We have a morphism of exact triangles
in $\lqcoh (Z)$ 
$$\xymatrix{\mathbb{T}_{Z} \ar[dd]_-{Q^{\flat}} \ar[r] & 
pr_{1}^{*}(\mathbb{T}_{Y_1}) \oplus pr_{2}^{*}(\mathbb{T}_{Y_2}) \ar[r] \ar[dd]_-{pr^*_1\Theta_{\gamma_1} \oplus
pr^*_2\Theta_{\gamma_2}} & p_1^{*}(\mathbb{T}_{X}) \ar[dd]^-{p_1^*(q^{\flat})} \\ & & \\
\mathbb{L}_{Z}[n-1] \ar[r] & pr_{1}^{*}(\mathbb{L}_{f_1})[n-1] \oplus pr_{2}^{*}(\mathbb{L}_{f_2})[n-1] \ar[r]
& p_{1}^*(\mathbb{L}_{X}[n]).}$$
Now, the morphism $p_1^*(q^{\flat})$ is a quasi-isomorphism since $q$ is non-degenerate, and  
the morphism $pr^*_1\Theta_{\gamma_1} \oplus
pr^*_2\Theta_{\gamma_2}$ is a quasi-isomorphism because the two null-structures are non-degenerate. This
implies that $Q^{\flat}$ is a quasi-isomorphism too, and thus that 
$Q=Q(q,\gamma_1,\gamma_2) \in \mathrm{QF}^{\textrm{nd}}(Z;n-1)$. 

\hfill $\Box$





\section{Derived Clifford algebra of derived quadratic complexes and stacks}

This section is rather brief and sketchy, since we mainly observe that the main definitions and basic results of Section \ref{dercliff} go through over a base derived Artin stack $X$. \\
All derived stacks even when not explicitly stated will be Artin and locally of finite presentation over $k$.\\

Let $X$ be a derived Artin stack locally finitely presented over $k$. We denote by $\mathsf{Alg}_{X}$ (respectively, $\mathsf{Alg}^{\textrm{perf}}_{X}$) the $\infty$-category of associative algebra objects in the symmetric monoidal $\infty$- category $(\lqcoh (X), \otimes_{{O}_{X}})$ (respectively, $(\lperf (X), \otimes_{{O}_{X}})$. Its objects will be simply called \emph{Algebras} (respectively, \emph{perfect Algebras}) over $X$. For any $\infty$-category $\mathsf{T}$, and any pair $(x,y)$ of objects in $\mathsf{T}$, we denote by $\mathrm{Map}_{\mathsf{T}}(x,y)$ the corresponding mapping space. Note that, up to isomorphisms in $\ho (\ssets)$, this is compatible withe notations used before when $\mathsf{T}$ is the Dwyer-Kan localization of a model category with respect to weak equivalences.\\

Exactly as in the case of $X$ being the derived spectrum of a simplicial commutative $k$-algebra, for \emph{evenly} shifted derived quadratic complexes over $X$ , it is possible to define a derived Clifford Algebra. We will sketch briefly the definitions and results, leaving to the reader the necessary changes with respect to the derived affine case.\\ Let $n\in \mathbb{Z}$, and $(E,q)$ be derived $2n$-shifted quadratic complex over $X$. 
The derived Clifford algebra functor associated to $(E,q)$ is defined by $$\underline{\mathbf{Cliff}}_{X}(E,q,2n): \mathsf{Alg}_{X} \longrightarrow \ssets \,\, : \, B \longmapsto \underline{\mathbf{Cliff}}(E,q,2n)(B)$$ where $\underline{\mathbf{Cliff}}(E,q,2n)(B)$ is defined by the following homotopy pull back in $\ssets$ $$\xymatrix{\underline{\mathbf{Cliff}}(E,q,2n)(B) \ar[r] \ar[d] & \mathrm{Map}_{\lqcoh (X)}(E,B[n]) \ar[d]^-{s_B} \\ \textrm{*} \ar[r]_-{\widetilde{q}_{B}} &  \mathrm{Map}_{\lqcoh (X)}(E\otimes^{}_{A}E,B[2n])}$$ where the maps $s_{B}$ and $q_{B}$ are defined analogously as in the case where $X= \mathbb{R}\mathsf{Spec} (A)$, for $A\in \salg$.

\begin{prop}\label{cliff} The functor $\underline{\mathbf{Cliff}}(E,q,2n)$ is co-representable, i.e. there exists a well defined $\mathbf{Cliff}_{A}(E,q,2n) \in \ho (\mathsf{Alg}_{X})$ and a canonical isomorphism in $\ho (\ssets)$ $$\underline{\mathbf{Cliff}}(E,q,2n)(B)\simeq \mathrm{Map}_{\mathsf{Alg}_{X}}(\mathbf{Cliff}_{A}(E,q,2n),B).$$
\end{prop}

\noindent \textbf{Proof.}   The proof is analogous to the one of Proposition \ref{Cliff}.
\hfill $\Box$ \\
 
 \begin{df} The Algebra $\mathbf{Cliff}_{X}(E,q;2n)$, defined up to isomorphism in $\mathsf{Alg}_{X}$, is called the \emph{derived Clifford Algebra} of the derived $2n$-shifted quadratic complex $(E,q)$. When $E=\mathbb{T}_{X}$, we will write $$\mathbf{Cliff}(X,q;2n):= \mathbf{Cliff}_{X}(\mathbb{T}_{X},q;2n).$$
 \end{df}
 
When the base derived stack $X$ is clear from the context, we will simply write $\mathbf{Cliff}(E,q,2n)$ for $\mathbf{Cliff}_{X}(E,q,2n)$.





\begin{prop}\label{presglobal}  Let $X$ be a derived Artin stack, lfp over $k$, $n \in \mathbb{Z}$, $f: E_1 \to E_2$ be a map in $\lqcoh (X)$, and $q_2$ a derived $2n$-shifted quadratic form on $E_2$ over $X$. Then there is a canonical map in $\ho (\dga)$ $$\mathbf{Cliff}_{X}(E_1,f^*q_2,2n) \longmapsto \mathbf{Cliff}_{X}(E_2,q_2,2n)$$ where $f^*q_2$ is the pull-back quadratic form of Definition \ref{resglobal}.
\end{prop}

\noindent \textbf{Proof.}   Same proof as for Proposition \ref{pres}.\hfill $\Box$ \\

Note that, by composition with the natural adjunction map in $\lqcoh (X)$ $$E[-n] \rightarrow \mathsf{Free}_{\mathcal{O}_{X}}(E[-n]),$$ the derived Clifford algebra $\mathbf{Cliff}_{X}(E,q,2n)$ of  the derived $2n$-shifted quadratic complex $(E,q)$, comes equipped with a natural map in $\lqcoh (X)$ $$E[-n] \longrightarrow \mathbf{Cliff}_{X}(E,q,2n).$$

The base-change of derived quadratic complexes over stacks (along maps of derived stacks), and derived isometric structures on a map of derived quadratic complexes (over a fixed derived stack), induce the following behavior on derived Clifford Algebras.

\begin{prop} Let $n \in \mathbb{Z}$. \begin{enumerate}
\item  If $(E_1,q_1)$ and $(E_2,q_2)$ are $2n$-shifted derived quadratic complexes over $X$, $f:E_1 \to E_2$ a map in $\lqcoh (X)$ and $\gamma \in  \pi_{0}(\mathsf{Isom}(f;(E_1,q_1),(E_2,q_2)))$ is a derived isometric structure on $f$, then there is an induced map in $\lqcoh (X)$ $$f_{\gamma}: \mathbf{Cliff}_{X}(E_1,q_1,2n) \longrightarrow \mathbf{Cliff}_{X}(E_2,q_2,2n).$$ If moreover $f$ is a quasi-isomorphism, $f_{\gamma}$ is an isomorphism.
\item If $(E,q)$ is a $2n$-shifted derived quadratic complex over $Y$, and $\varphi :X \to Y$ a morphism in $\dst$. Then there is a canonical isomorphism in $\lqcoh (X)$ $$\xymatrix{\varphi^{*} (\mathbf{Cliff}_{Y}(E,q,2n)) \ar[r]^{\sim} & \mathbf{Cliff}_{X}(\varphi^{*}E,\varphi^{*}q,2n)}.$$ 
\item Let $\varphi: X \to Y$ be a map in $\dst$, and $q$ a derived $2n$-shifted quadratic form on $Y$. By Definition \ref{shriek},  $(X, \varphi^{q} q)$ is a derived $2n$-shifted quadratic stack, and there is a canonical map in $\mathsf{Alg}_{X}$ $$\mathbf{Cliff}(X,\varphi^{q} q;2n) \longrightarrow \varphi^{*}\mathbf{Cliff}(Y,q;2n)$$
\end{enumerate}
\end{prop}

\noindent \textbf{Proof.}  (1) and (2) are analogous to the proofs of propositions \ref{p1} and \ref{p2}, respectively. To establish the map in (3), we first observe that the map $\mathbb{T}_{X} \to \varphi^* \mathbb{T}_{Y}$ gives us, by Prop. \ref{presglobal}, an induced canonical map 
$$\mathbf{Cliff}(X,\varphi^{q} q;2n) =\mathbf{Cliff}_{X}(\mathbb{T}_{X},\varphi^{q} q;2n) \longrightarrow \mathbf{Cliff}_{X}(\varphi^* \mathbb{T}_{Y},\varphi^{*} q;2n)$$ and we conclude by point (2) applied to $E= \mathbb{T}_{Y}$.
 \hfill $\Box$ \\

\begin{rmk} \emph{In the situation, and notations, of Theorem \ref{map}, with both $n$ and $d$ even, it can be proved that there exists a canonical map in $\mathsf{Alg}_{\textsf{MAP}_{\dst}(Y,X)}$ $$h_{\eta}: \mathbf{Cliff}(\textsf{MAP}_{\dst}(Y,X) ,q'_{\eta}; n-d) \longrightarrow \mathrm{pr}_{*}\mathrm{ev}^* \mathbf{Cliff}(X, q;n)$$ where $\mathrm{pr}: Y \times \textsf{MAP}_{\dst}(Y,X) \to \textsf{MAP}_{\dst}(Y,X)$ is the projection, and $\mathrm{ev}: Y \times \textsf{MAP}_{\dst}(Y,X) \to X$  the evaluation map. It is an open question whether this map is an isomorphism in $\mathsf{Alg}_{\textsf{MAP}_{\dst}(Y,X)}$. }
\end{rmk}

\section{Appendix: Superstuff}\label{app}

In this Appendix, we recall a few basic facts about $\mathbb{Z}$- and $\mathbb{Z}/2$-graded dg-modules and dg-algebras, mainly to fix our notations and establish the background for Proposition \ref{graded} in the main text, where we need to know that the derived Clifford algebra is naturally an object in the homotopy category of $\mathbb{Z}/2$-graded dg-algebras. In order to distinguish the internal grading from the external one, we will call the former the (cohomological) \emph{degree}, and the latter the \emph{weight}.\\
 
Let $A \in \cdga$, and $A-\dgmod^{w}\equiv A-\dgmod^{\mathbb{Z}/2-\textrm{gr}}$ denote the category of $\mathbb{Z}/2$\emph{-weighted $A$-dg modules}, whose objects are triples $C_{*}^{\bullet}=(C; C_0,C_1)$ where $C \in A-\dgmod$, and $(C_0,C_1)$ are sub $A$-dg modules that provide a $\mathbb{Z}/2$-graded decomposition in $A-\dgmod$ $C=C_0\oplus C_1$. Elements of $C_0$ (resp. of $C_1$) are said to have \emph{weight} $0$ (respectively, $1$). The morphisms in  $A-\dgmod^{w}$ $(C; C_0, C_1)\to (D; D_0, D_1)$ are morphisms $f:C\to D$ in $A-\dgmod$ preserving the weights: $f(C_i)\subset D_i$, $i=0,1$.
There is a symmetric monoidal structure on $A-\dgmod^{w}$, defined by $$(C; C_0, C_1) \otimes_{A}^{w} (C; C_0, C_1):= (C\otimes_{A} D, (C_0 \otimes_{A} D_0) \oplus (C_1 \otimes_{A} D_1), (C_0 \otimes_{A} D_1) \oplus (C_1 \otimes_{A} D_0) ),$$ with commutativity constraint given, on $(\mathbb{Z} \times \mathbb{Z}/2)$-homogeneous elements, by $$\sigma_{(C_{*}^{\bullet},D_{*}^{\bullet})}: C_{*}^{\bullet} \otimes^{w}D_{*}^{\bullet} \longrightarrow  D_{*}^{\bullet} \otimes^{w}C_{*}^{\bullet}\, : \, x\otimes y \longmapsto (-1)^{w(x)w(y) + \deg (x) \deg (y)} y\otimes x$$ where $w(-)$ denotes the weight, and $\deg(-)$ the (cohomological) degree. Alternatively, we could have  written $$\sigma_{(C_{*}^{\bullet},D_{*}^{\bullet})} (x_{w}^i \otimes y_{w'}^{j})= (-1)^{ww'+ ij} y_{w'}^{j} \otimes x_{w}^{i},$$ with the standard obvious meaning of the symbols.\\
There are functors $$\mathsf{i}: A-\dgmod \longrightarrow A-\dgmod^{w} \, : C \longmapsto (C; C_0=C, 0)\, ,$$ $$\textsf{forget}: A-\dgmod^{w} \longrightarrow A-\dgmod \, : (C,C_0,C_1) \longmapsto C\, , $$ but observe that the first one is symmetric monoidal, while the second one is not.\\

The category $A-\dgmod^{w}$ is endowed with a symmetric monoidal model structure where equivalences and fibrations are morphisms mapped to quasi-isomorphisms and fibrations via the forgetful functor to $A-\dgmod$: the classical proof in the unweighted case goes through. This model structure satisfies the monoidal axiom (\cite[Def. 3.3]{ss1}), hence the category $A-\dga^{w}$ of monoids in $(A-\dgmod^{w}, \otimes^{w})$ have an induced model category structure where equivalences and fibrations are detected on the underlying $\mathbb{Z}/2$-weighted $A$-dg modules (\cite[Thm. 4.1]{ss1}). The objects of $A-\dga^{w}$ will be called $\mathbb{Z}/2$\emph{-weighted $A$-dg algebras}. Note that both $A-\dgmod^{w}$ and $A-\dga^{w}$ are cofibrantly generated model categories. The forgetful functor $$\textsf{Forget}: A-\dga^{w} \longrightarrow A-\dgmod^{w}$$ is Quillen right adjoint to the free $\mathbb{Z}/2$-weighted $A$-dg algebra functor $$\mathsf{T}_{A}: A-\dgmod^{w} \longrightarrow A-\dga^{w}$$ (given by the tensor algebra construction in $(A-\dgmod^{w}, \otimes^{w})$). Note that $T$ acquires also an additional $\mathbb{Z}$-grading but we will not use it.\\ Note that moreover, since $A$ is (graded commutative), the category $\ho (A-\dga^{w})$ comes equipped with a (derived) tensor product, defined by $$(B;B_0,B_1)\otimes_{A}^{w} (D;D_0,D_1) := Q(B;B_0,B_1)\otimes_{A}^{w} (D;D_0,D_1)$$ where $Q(-)$ denotes a cofibrant replacement functor in $A-\dga^{w}$, and the algebra product is defined by $$(x_w^i \otimes y_{w'}^j)\cdot (z_{p}^h\otimes t_{q}^{k}):= (-1)^{w'p}(x_w^i z_{p}^h)\otimes (y_{w'}^j t_{q}^{k}).$$
\linebreak

Completely analogous (notations and) results hold if we start with the category $A-\dgmod^{\mathbb{Z}-\textrm{gr}}$ of $\mathbb{Z}$\emph{-weighted $A$-dg modules}, whose objects are pairs $C_{*}^{\bullet}=(C; C_{w})_{w \in \mathbb{Z}}$ where $C \in A-\dgmod$, and $(C_{w})_{w \in \mathbb{Z}}$ are sub $A$-dg modules providing a $\mathbb{Z}$-graded decomposition in $A-\dgmod$ $C=\oplus_{w}C_w$. The symmetric monoidal structure is given by $$(C; C_{w})_{w \in \mathbb{Z}} \otimes_{A}^{w} (D; D_{w'})_{w' \in \mathbb{Z}}:= (C\otimes_{A} D, ((C\otimes_{A} D)_{p}:=\oplus_{w+w'=p}C_{w} \otimes_{A} D_{w'})_{p \in \mathbb{Z}}) $$ while the commutativity is given by the same formula as in the $\mathbb{Z}/2$-graded case (except that the weights are now in $\mathbb{Z}$).
Thus we dispose of a (cofibrantly generated) model category $A-\dga^{\mathbb{Z}-\textrm{gr}}$ of monoids in the symmetric monoidal (cofibrantly generated) model category $(A-\dgmod^{\mathbb{Z}-\textrm{gr}}, \otimes^{w})$.\\

\begin{prop} The functor $\mathsf{Free}_{A}: A-\dgmod \rightarrow A-\dga$ (used in the main text) naturally factors through the weight-forgetful functor $A-\dga^{\mathbb{Z}-\textrm{gr}} \longrightarrow A-\dga$ (that simply forgets the weight-grading). And, 
if we denote by the same symbol the induced functor $$\mathsf{Free}_{A}: A-\dgmod \rightarrow A-\dga^{\mathbb{Z}-\textrm{gr}},$$  this functor sends weak equivalences between cofibrant $A$-dg modules to weak equivalences.
\end{prop}

\noindent \textbf{Proof.} Let $C \in A-\dgmod$. By giving to $\mathsf{Free}_{A}(C)$ the $\mathbb{Z}$-grading $$(\mathsf{Free}_{A}(C))_{w}:=\mathsf{T}_{A}^w(C)= C\otimes_{A} \cdots \otimes_{A}C \,\,\, \textrm{($w$ times)}$$ with the convention that $\mathsf{T}_{A}^w(C)=0$ for $w<0$, it is easy to verify that $$(\mathsf{Free}_{A}(C), (\mathsf{Free}_{A}(C))_{w})_{w \in \mathbb{Z}})$$ defines an element in $A-\dga^{\mathbb{Z}-\textrm{gr}}$. The induced functor $\mathsf{Free}_{A}: A-\dgmod \rightarrow A-\dga^{\mathbb{Z}-\textrm{gr}}$ sends weak equivalences between cofibrant objects to weak equivalences since the left Quillen functor $\mathsf{Free}_{A}: A-\dgmod \rightarrow A-\dga$ does (\cite[Lemma 1.1.12]{ho}).
\hfill $\Box$\\

By definition of $\otimes^{w}$ both in the $\mathbb{Z}/2$- and in the $\mathbb{Z}$-weighted case, the functor $$A-\dgmod^{\mathbb{Z}-\textrm{gr}} \longrightarrow A-\dgmod^{w} $$ given by $$(C; C_{w})_{w\in \mathbb{Z}} \longmapsto (C; C_\mathrm{even}:= \oplus_{w} C_{2w}, C_{\mathrm{odd}}:= \oplus _{w} C_{2w +1})$$ is symmetric monoidal, and preserves weak equivalences, hence it induces a weak-equivalences preserving functor between the corresponding monoid objects 
$$(-)^{w}: A-\dga^{\mathbb{Z}-\textrm{gr}} \longrightarrow A-\dga^{w}.$$
Thus the composite functor $$\xymatrix{\mathsf{Free}^{w}_{A}: A-\dgmod \ar[r]^-{\mathsf{Free}_{A}} & A-\dga^{\mathbb{Z}-\textrm{gr}} \ar[r]^-{(-)^{w}} &  A-\dga^{w}}$$ sends weak equivalences between cofibrant objects to weak equivalences. \\

Now we come back to explain the statement of Proposition \ref{graded}. Let $(C,q)$ be a derived $2n$-shifted quadratic complex over $A$, where $C$ is cofibrant in $A-\dgmod$. Note that, then, also any shift of $C$, and any shift of $C\otimes_A C$ is cofibrant, too.  Consider the homotopy push-out square in $A-\dga$ $$\xymatrix{\mathsf{Free}_{A}(C\otimes^{}_{A}C[-2n]) \ar[r]^-{u} \ar[d]_-{t} & \mathsf{Free}_{A}(C[-n]) \ar[d] \\ \mathsf{Free}_{A}(0)=A \ar[r] & \mathbf{Cliff}(C,q; 2n) }$$
used in the proof of Proposition \ref{cliff} in order to define the derived Clifford algebra $\mathbf{Cliff}(C,q; 2n)$.\\  By upgrading $\mathsf{Free}_{A}$ to $\mathsf{Free}^{w}_{A}$, we observe that the maps $u$ and $t$ are no longer maps in $A-\dga^{w}$ (they have odd degree) but they both land into the even parts of $\mathsf{Free}^{w}_{A}(C[-n])$ and $\mathsf{Free}^{w}_{A}(0)$, respectively. This implies that $\mathbf{Cliff}(C,q; 2n)$ admits a $\mathbb{Z}/2$-weight grading $$(\mathbf{Cliff}(C,q; 2n); \mathbf{Cliff}(C,q; 2n)_0, \mathbf{Cliff}(C,q; 2n))$$ such that the induced maps $$\mathsf{Free}^{w}_{A}(C[-n]) \longrightarrow (\mathbf{Cliff}(C,q; 2n); \mathbf{Cliff}(C,q; 2n)_0, \mathbf{Cliff}(C,q; 2n)),  $$ $$\mathsf{Free}^{w}_{A}(A)=(A;A,0) \longrightarrow (\mathbf{Cliff}(C,q; 2n); \mathbf{Cliff}(C,q; 2n)_0, \mathbf{Cliff}(C,q; 2n)),$$ are maps in $\ho (A-\dga^{w})$. In other words the derived Clifford algebra functor $$(C,q:2n)\longmapsto \mathbf{Cliff}(C,q; 2n)$$ extend to a functor $$\mathbf{Cliff}^{w}_{A}: \underline{\mathsf{QMod}}(A; 2n) \longrightarrow  L^{DK}(A-\dga^{w}),$$ where 
\begin{itemize}
\item $\underline{\mathsf{QMod}}(A; 2n)$ is the $\infty$-category of Def. \ref{qmodA}, and
\item $L^{DK}(A-\dga^{w})$ is the Dwyer-Kan localization of the model category $A-\dga^{w}$ along its weak equivalences.
\end{itemize}


\begin{thebibliography}{50}

\bibitem[Be-Fa]{befa}K. Behrend, B. Fantechi, \emph{Symmetric obstruction theories and Hilbert schemes of points on threefolds},
Algebra Number Theory, 2 (2008), 313–345.

\bibitem[Do]{dold} A. Dold, \emph{Homology of symmetric products and other functors of complexes}, Ann. Math. 68 (1) (1958) 54–80

\bibitem[D-K]{dk} W. Dwyer, D. Kan, \emph{Simplicial localizations of categories}, J. Pure Appl. Algebra 17 (1980), 267–284

\bibitem[Fr-Ha]{fricke} J. Fricke, L. Habermann, \emph{
On the geometry of moduli spaces of symplectic structures}, Manuscripta Math.
December 2002, Volume 109, Issue 4, pp 405-417.

\bibitem[G-K]{gaka} N. Ganter, M. Kapranov \emph{Symmetric and exterior powers of categories}, Preprint arXiv:1110.4753v1.

\bibitem[Ho]{ho} M. Hovey, \textit{Model categories}, Mathematical surveys and monographs, Vol. $\mathbf{63}$,
Amer. Math. Soc., Providence 1998.

\bibitem[Lu--\textsf{HTT}]{htt} J. Lurie, Higher Topos Theory, AM - 170, Princeton University Press, Princeton, 2009

\bibitem[Ma]{malm} E. J. Malm, \emph{String topology and the based loop space}, Preprint arXiv:1103.6198.

\bibitem[Mi-Re]{mire} A. Micali, Ph. Revoy, \emph{Modules quadratiques}, Mémoires de la SMF, tome 63 (1979), p. 5-144.
 



\bibitem[PTVV]{ptvv}  T. Pantev, B. To\"en, M. Vaqui\'e, G. Vezzosi, \emph{Shifted Symplectic Structures},  Publications Math\'ematiques de l'IHES, 
June 2013, Volume 117, Issue 1, pp 271-328, DOI: 10.1007/s10240-013-0054-1

\bibitem[Pre]{pre} A. Preygel, \emph{Thom-Sebastiani \& Duality for Matrix Factorizations}, Preprint arXiv:1101.5834.

\bibitem[Ra]{ra} A. Ranicki, \emph{Algebraic L-theory and Topological Manifolds}, Cambridge Tracts in Mathematics 102, CUP (1992).

\bibitem[Schl]{sch} M. Schlichting, \emph{Hermitian K-theory, derived equivalences, and Karoubi's Fundamental Theorem, } Preprint arXiv:1209.0848.

\bibitem[Sch-Shi-1]{ss1}S. Schwede, B. Shipley, \emph{Algebras and modules in monoidal model categories}, Proc. London Math. Soc. (3)
80 (2000), no. 2, 491-511.

\bibitem[Sch-Shi-2]{ss2} S. Schwede, B. Shipley, \emph{Equivalences of monoidal model categories}, Algebr. Geom. Topol. 3 (2003), 287-
334.

\bibitem[Sim]{simp} C. Simpson, \emph{Geometricity of the Hodge filtration on the 1-stack of perfect complexes
over $X_{DR}$}, Moscow Mathematical Journal, Vol. 9 (2009), 665-721

\bibitem[Su]{su} D.Sullivan, \emph{Infinitesimal computations in topology}, Publ. Math. I.H.E.S.  47 (1977) pp. 269-331.

\bibitem[Th]{tt} T. Thomas, \emph{Cobordism, categories, quadratic forms}, talk transcript available at \texttt{http://www.maths.ed.ac.uk/$\sim$aar/papers/toptalk.pdf}.

\bibitem[To--\textsf{Seattle}]{seattle} B. To\"en, \textit{Higher and derived stacks: a global overview}, Algebraic Geometry—Seattle 2005. Part 1, Proc.
Symp. Pure Math., vol. 80, pp. 435–487, Am. Math. Soc., Providence, 2009.

\bibitem[To--\textsf{deraz}]{deraz} B. To\"en, \emph{Derived Azumaya algebras and generators for twisted derived categories}, Invent. Math. 189 (2012), no. 3, 581–652.






\bibitem[To-Ve --\textsf{Traces}]{chern} B. To\"en, G. Vezzosi, \emph{Caractères de Chern, traces équivariantes et géométrie
algébrique dérivée},  to appear in Selecta Math.

\bibitem[To-Ve --\textsf{HAGII}]{hagII} B. To\"en, G. Vezzosi, \textit{Homotopical algebraic geometry II: Geometric stacks
and applications},  
Mem. Amer. Math. Soc.  \textbf{193}  (2008),  no. 902, x+224 pp.


\end{thebibliography}
\end{document}